\documentclass[toc=bibnumbered,ngerman,english,titlepage=false]{scrartcl}

\usepackage[utf8]{inputenc}
\usepackage{babel}
\usepackage[T1]{fontenc}
\usepackage{lmodern}
\usepackage{microtype}

\usepackage{amsmath}
\usepackage{amsfonts}
\usepackage{amssymb}
\usepackage{bm}
\usepackage{textcomp}
\usepackage{csquotes}
\usepackage{enumitem}
\usepackage{graphicx}

\usepackage{scrpage2}

\numberwithin{equation}{section} 

\newcommand{\tildeTheta}{\tilde{\Theta}}

\newcommand{\N}{\ensuremath{\mathbb{N}}}

\newcommand{\R}{\ensuremath{\mathbb{R}}}

\newcommand{\C}{\ensuremath{\mathbb{C}}}
\newcommand{\K}{\ensuremath{\mathbb{K}}}

\newcommand{\E}{\ensuremath{\mathbb{E}}}
\newcommand{\B}{\ensuremath{\mathbb{B}}}
\newcommand{\diverg}{\textup{div}}

\newcommand{\tr}{\textup{tr}}
\newcommand{\Tr}{\textup{Tr}}

\newcommand{\lsprung}{\textup{\textlbrackdbl}}
\newcommand{\rsprung}{\textup{\textrbrackdbl}}

\newcommand{\Lc}{\ensuremath{\mathcal{L}}}

\setkomafont{disposition}{\normalcolor\bfseries}
 \setkomafont{section}{\large} 
 
\begin{document}
\pagenumbering{roman}

\oddsidemargin 6pt\evensidemargin 6pt\marginparwidth 48pt\marginparsep 10pt

\topmargin -18pt\headheight 12pt\headsep 25pt

\ifx\cs\documentclass \footheight 12pt \fi \footskip 30pt

\textheight 625pt\textwidth 431pt\columnsep 10pt\columnseprule 0pt 
                                                                      
\automark[section]{section}
\clearscrheadings
\chead{\scshape\headmark}
\cfoot{\pagemark}
\pagestyle{scrheadings}

 \renewcommand{\headfont}{\slshape}      
 \renewcommand{\pnumfont}{\upshape}      
 \setcounter{secnumdepth}{4}             
 \setcounter{tocdepth}{4}             

\newtheorem{Definition}{Definition}[section]
\newtheorem{Satz}[Definition]{Satz}
\newtheorem{Lemma}[Definition]{Lemma}
\newtheorem{Korollar}[Definition]{Korollar}
\newtheorem{Corollary}[Definition]{Corollary}
\newtheorem{Bemerkung}[Definition]{Bemerkung}
\newtheorem{Remark}[Definition]{Remark}
\newtheorem{Proposition}[Definition]{Proposition}
\newtheorem{Beispiel}[Definition]{Beispiel}
\newtheorem{Theorem}[Definition]{Theorem}

\begin{center}
\textbf{\Large Well-Posedness of a Navier-Stokes/Mean Curvature Flow system}\\
\vspace{0.5cm}
\large\textsc{Helmut Abels\quad and\quad Maximilian Moser}\\
\vspace{0.2cm}
\textit{Fakultät für Mathematik, Universität Regensburg, Universitätsstraße 31, D-93053 Regensburg, Germany}\\
\vspace{0.2cm}
\small E-mail: \textsf{Helmut.Abels@mathematik.uni-regensburg.de}, 
\textsf{Maximilian1.Moser@mathematik.uni-regensburg.de}\\
\vspace{0.2cm}
\normalsize \today
\end{center}
\textbf{Abstract.} We consider a two-phase flow of two incompressible, viscous and immiscible fluids which are separated by a sharp interface in the case of a simple phase transition. In this model the interface is no longer material and its evolution is governed by a convective mean curvature flow equation, which is coupled to a two-phase Navier-Stokes equation with Young-Laplace law. The problem arises as a sharp interface limit of a diffuse interface model, which consists of a Navier-Stokes system coupled with an Allen-Cahn equation. We prove existence of strong solutions for sufficiently small times and regular initial data.

\vspace{0.1in}\\

 \vspace{0.1in}\\

\pagenumbering{arabic}

\setcounter{secnumdepth}{4}\setcounter{tocdepth}{4} 
\section{Introduction}\label{sec_einl}\thispagestyle{plain}
We consider the flow of two immiscible, incompressible Newtonian fluids with phase transition in a bounded, smooth domain $\Omega\subseteq\R^n,n=2,3$. At time $t\in(0,T)$ the fluids fill domains $\Omega^\pm(t)$ separated by the interface $\Gamma(t):=\partial\Omega^+(t)$. For simplicity we set the densities to one. Moreover, we assume that the fluids have constant viscosities $\mu^\pm>0$. Then the stress tensor has the form $T(v,p)=2\mu^\pm \textup{sym}(\nabla v)-pI$, where $\textup{sym}(\nabla v):=\frac{1}{2}(\nabla v+\nabla v^\top)$, $v:\Omega\times(0,T)\rightarrow\R^n$ is the velocity and $p:\Omega\times(0,T)\rightarrow\R$ is the pressure in Eulerian coordinates.
 
To formulate the model, we need some notation. We denote by $\nu_{\Gamma(t)}$ the outer unit normal to $\Gamma(t)=\partial\Omega^+(t)$ and by $V_{\Gamma(t)}, H_{\Gamma(t)}$ the normal velocity and mean curvature (for convenience the sum of the principal curvatures) with respect to $\nu_{\Gamma(t)}$. Furthermore the jump of a quantity $f$ defined on $\Omega^\pm(t)$ across 
$\Gamma(t)$ (with respect to $\nu_{\Gamma(t)}$) is defined as
\[
\lsprung f\rsprung(x):=\textup{lim}_{r\rightarrow0}f(x+r\nu_{\Gamma(t)}(x))-f(x-r\nu_{\Gamma(t)}(x))
\]
for all $x\in\Gamma(t)$.

As model we consider the following Navier-Stokes/mean curvature flow system
\begin{alignat}{2}
\partial_t v+v\cdot\nabla v-\mu^\pm\Delta v+\nabla p&=0&\qquad &\text{ in }\Omega^\pm(t),t\in(0,T),\label{eq_nsmc_1}\\
\diverg\,v&=0&\qquad &\text{ in }\Omega^\pm(t),t\in(0,T),\label{eq_nsmc_2}\\
-\lsprung T(v,p)\rsprung\,\nu_{\Gamma(t)}&=\sigma H_{\Gamma(t)}\nu_{\Gamma(t)}  &\qquad&\textup{ on }\Gamma(t),t\in(0,T),\label{eq_nsmc_3}\\
\lsprung v\rsprung&=0&\qquad  &\textup{ on }\Gamma(t), t\in(0,T),\label{eq_nsmc_4}\\
v|_{\partial\Omega}&=0 &\qquad &\textup{ on }\partial\Omega\times(0,T),\label{eq_nsmc_5}\\
v|_{t=0}&=v_0 &\qquad&\text{ in }\Omega,\label{eq_nsmc_6}\\
V_{\Gamma(t)}-\nu_{\Gamma(t)}\cdot v|_{\Gamma(t)}&=H_{\Gamma(t)} &\qquad &\textup{ on }\Gamma(t),t\in(0,T),\label{eq_nsmc_7}\\
\Gamma(0)&=\Gamma_0\label{eq_nsmc_8},&\qquad &
\end{alignat}
where $\sigma>0$ is a surface tension constant and $\Gamma_0$ is the initial interface. Here the velocity $v$, the pressure $p$ and the family of interfaces $\{\Gamma(t)\}_{t\in(0,T)}$ are to be determined. Equations \eqref{eq_nsmc_1}-\eqref{eq_nsmc_2} describe the (local) mass and momentum conservation of the fluids and \eqref{eq_nsmc_3} is the balance of forces at the interface, see Landau and Lifshitz \cite{Landau_Lifshitz_Fluid}, §61. Moreover we require continuity of $v$ on $\Gamma(t)$ in \eqref{eq_nsmc_4} and the no-slip condition for $v$ on $\partial\Omega$ in \eqref{eq_nsmc_5}. \eqref{eq_nsmc_6} and \eqref{eq_nsmc_8} are the initial conditions for $v$ and the family of hypersurfaces $\{\Gamma(t)\}_{t\in(0,T)}$, respectively. Equation \eqref{eq_nsmc_7} without the curvature term $H_{\Gamma(t)}$ means that the interface is just transported by the fluids. Then \eqref{eq_nsmc_1}-\eqref{eq_nsmc_8} reduces to the classical two-phase Navier-Stokes system with surface tension, which was studied e.g. by Köhne, Prüss and Wilke \cite{Koehne_Pruess_Wilke_2N}. For this system the total mass of each fluid is preserved since the densities are constant and the enclosed volume of $\Gamma(t)$ stays the same. Formally this follows from (references in) Escher and Simonett \cite{EscherSimonett_HeleShaw}, p. 623, the divergence theorem and \eqref{eq_nsmc_2}:
\[
\frac{d}{dt}|\Omega^+(t)|=\int_{\Gamma(t)}V_{\Gamma(t)}\,d\sigma=\int_{\Gamma(t)}v\cdot\nu_{\Gamma(t)}\,d\sigma=\int_{\Omega^+(t)}\diverg\,v\,dx=0.
\]
The additional curvature term $H_{\Gamma(t)}$ on the right hand side in \eqref{eq_nsmc_7} allows for a change of total masses and couples the Navier-Stokes system to the mean curvature flow
\[
V_{\Gamma(t)}=H_{\Gamma(t)}\quad\text{ for }t>0.
\]
Therefore we call \eqref{eq_nsmc_1}-\eqref{eq_nsmc_8} the Navier-Stokes/mean curvature flow system. The system is obtained as a sharp interface limit of a Navier-Stokes-Allen-Cahn system. Another motivation is the regularizing effect of the curvature term in \eqref{eq_nsmc_7}. But we can also view it as a simple model for a two-phase flow with phase transition. In the case that the Navier-Stokes equations are modified for shear thickening non-Newtonian fluid of power-law type  Liu, Sato, and Tonegawa in \cite{LiuSatoTonegawa2} have shown existence of weak solutions. In the present contribution we prove local existence of strong solutions to \eqref{eq_nsmc_1}-\eqref{eq_nsmc_8} in an $L^p$-setting. The result is part of the second author's Master thesis, which was supervised by the first author. 

We apply a similar strategy as A. and Wilke \cite{Abels_Mullins_Sekerka}, who study local well-posedness and qualitative behaviour of solutions for a two-phase Navier-Stokes-Mullins-Sekerka system. The equations are transformed with a variant $\tilde{\Theta}_h$ of the Hanzawa transform to a reference surface $\Sigma$ and disjoint domains $\Omega^\pm$, where $\Omega=\Omega^+\cup\Omega^-\cup\Sigma$ is similar as above. Here $h:\Sigma\times(0,T)\rightarrow\R$ is the \enquote{height function} of $\{\Gamma(t)\}_{t\in(0,T)}$ with respect to $\Sigma$. The coupling with the Navier-Stokes part in \eqref{eq_nsmc_7} is of lower order. Therefore we solve the Navier-Stokes-part \eqref{eq_nsmc_1}-\eqref{eq_nsmc_6} for a given appropriate family of hypersurfaces $\{\Gamma(t)\}_{t\in(0,T)}$ or equivalently for a given suitable time-dependent height function $h$. Then we insert the obtained velocity field $u(h)$ into equation \eqref{eq_nsmc_7} and \eqref{eq_nsmc_7}-\eqref{eq_nsmc_8} turns into an abstract evolution equation for $h$. The latter is solved by using the theory of maximal $L^p$-regularity.

The paper is organized as follows. In Section \ref{sec_prelim} we describe the transformation to a fixed reference hypersurface (and related domains) and analyse the occurring terms, especially the transformed mean curvature. In Section \ref{sec_NavSt_teil} we consider the Navier-Stokes part for a given family of hypersurfaces $\{\Gamma_h(t)\}_{t\in(0,T)}$ and in Section \ref{sec_local_wellposed} we show the local well-posedness for the transformed system. Definitions and some properties of the used Banach-space valued function spaces are summarized in Appendix \ref{sec_fctsp_vecvalued}. In Appendix \ref{sec_lin_stokes} we prove a maximal regularity result for a two-phase Stokes system which is needed for Section \ref{sec_NavSt_teil}.

We always assume $\Omega\subseteq\R^n, n=2,3$, to be a bounded, connected and smooth domain and throughout let $\Sigma\subseteq\Omega$ be a compact, connected and smooth hypersurface, that separates $\Omega$ in two disjoint, connected domains $\Omega^\pm$ with $\Sigma=\partial\Omega^+$ and outer unit normal $\nu_\Sigma$.

\noindent
\textbf{Acknowledgements:} The second author is grateful to Mathias Wilke for many helpful discussions and his lecture series on ``maximal regularity''.

\section{Preliminaries}\label{sec_prelim}
\subsection{Notation and Scalar-Valued Function Spaces}
Let $\R_+:=[0,\infty)$ and $\R^n_+:=\R^{n-1}\times(0,\infty)$ for $n\in\N, n\geq 2$. The Euclidean norm in $\R^n$ is denoted by $|.|$. For matrices $A,B\in\R^{n\times n}$ let $A:B$ be the matrix-product and $|A|:=\sqrt{A:A}$ be the induced norm. The set of invertible matrices in $\R^{n\times n}$ is $GL(n,\R)$. For $A\in GL(n,\R)$ we set $A^{-\top}:=(A^{-1})^\top$.

For metric spaces $X$ denote by $B_X(x,r)$ the ball with radius $r>0$ around $x\in X$. For normed spaces $X,Y$ over $\K=\R$ or $\C$ the set of bounded, linear operators $T:X\rightarrow Y$ is $\Lc(X,Y)$ and $\Lc(X):=\Lc(X,X)$. If $X,Y$ are additionally complete and $U\subseteq X$ is open, for $k\in\N$ the set of all $k$-times continuously differentiable functions $f:U\rightarrow Y$ is denoted by $C^k(U;Y)$. The set of all bounded $f\in C^k(U;Y)$ with bounded derivatives is $BC^k(U;Y)$ with usual norm.

If $(X_0,X_1)$ is admissible and $\theta\in(0,1),1\leq p \leq \infty$, we denote by $(X_0,X_1)_{\theta,p}$ the real interpolation spaces and with $\|.\|_{(X_0,X_1)_{\theta,p}}\equiv\|.\|_{\theta,p}$ the norm, see Lunardi \cite{Lunardi_Interpolation}, Bergh and Löfström \cite{BL_Interpolation} and Triebel \cite{Triebel_Interpol_Theory}.

Let $\Omega\subseteq\R^n,n\in\N$, be open and nonempty. Spaces of continuous and continuously differentiable functions $f:\overline{\Omega}\rightarrow\R$ are defined as usual. Now let $1\leq p\leq\infty$ and $k\in\N$. Then the Lebesgue and Sobolev spaces are denoted by $L^p(\Omega)=:W_p^0(\Omega)$ and $W^k_p(\Omega)$, respectively. Moreover, we set
\[
W^k_{p,0}(\Omega):=\overline{C_0^\infty(\Omega)}^{W^k_p(\Omega)}\quad\text{ and }\quad W^{-k}_p(\Omega):=W^k_{p',0}(\Omega)',
\]
where $p'$ is the conjugate exponent to $p$, as well as $L^p_{(0)}(\Omega):=\{f\in L^p(\Omega):\int_{\Omega}f\,dx=0\}$,
\[
W^k_{p,(0)}(\Omega):=W^k_p(\Omega)\cap L^p_{(0)}(\Omega)\quad\text{ and }\quad W^{-k}_{p,(0)}(\Omega):=W^k_{p',(0)}(\Omega)'.
\]
The closure of the divergence-free $\phi\in C_0^\infty(\Omega)^n$ in $L^p(\Omega)^n$ is denoted by $L^p_\sigma(\Omega)$. Besides, we need the Besov spaces $B^s_{p,q}(\R^n)$ for $0<s<\infty$ and $1\leq p,q\leq\infty$. We set $W^s_p(\R^n):=B^s_{p,p}(\R^n)$ for all $s\in(0,\infty)\textbackslash\N,1\leq p\leq \infty$. Respective spaces can also be defined on bounded, smooth domains $\Omega\subseteq\R^n$ and, if $n\geq 2$, on their boundaries $\partial\Omega$. For embeddings, interpolation results and trace theorems, compare e.g. Triebel \cite{Triebel_Interpol_Theory}, \cite{Triebel_Fct_SpacesI}.

\subsection{Hanzawa Transformation $\Theta_h$}\label{sec_hanzawa}
For the construction of $\Theta_h$ we use that there is an $a>0$, such that 
\[
X:\Sigma\times(-a,a)\rightarrow\R^n:(s,r)\mapsto s+r\,\nu_\Sigma(s)
\]
is a $C^\infty$-diffeomorphism onto $\Sigma_a:=B_a(\Sigma)$. We set $(\Pi,d_\Sigma):=X^{-1}$. Then it holds $\nabla d_\Sigma=\nu_\Sigma\circ\Pi$. Here $\Pi$ is called projection onto $\Sigma$ and $d_\Sigma$ signed distance function. For the proof cf. Hildebrandt \cite{Hildebrandt_Analysis2}, Chapter 4.6, Theorems 1-3.

Let $a<\textup{dist}(\Sigma,\partial\Omega)$ and $a_0<\frac{a}{4}$ fixed, $\chi\in C^\infty(\R)$ be a cutoff function such that $|\chi'|\leq 4$ and $\chi(s)=1$ for $|s|\leq \frac{1}{3}$, as well as $\chi(s)=0$ for $|s|\geq \frac{2}{3}$. For $h\in C^2(\Sigma)$ with $\|h\|_\infty<a_0$ we define the \textit{Hanzawa transformation} $\Theta_h$ by
\[
\Theta_h(x):=x+\chi(\frac{d_\Sigma(x)}{a})\,h(\Pi(x))\,\nu_\Sigma(\Pi(x))\quad\text{ for all }x\in\R^n.
\]

Essential properties of $\Theta_h$ are listed in
\begin{Lemma}\label{th_hanzawa}
The Hanzawa transformation $\Theta_h:\Omega\rightarrow\Omega$ is a $C^2$-diffeomorphism and $\Theta_h$ is the identity on $\R^n\textbackslash \Sigma_{2a/3}$, in particular nearby $\partial\Omega$. Moreover,
\[
\det D\Theta_h\geq c>0\quad\text{ and }\quad\|D\Theta_h\|_{\infty}+\|D\Theta_h^{-1}\|_\infty\leq C(1+\|h\|_{C^1(\Sigma)})
\]
with $c,C>0$ independent of $h$.
\end{Lemma}
\textit{Proof.} The first part can be shown as in Kneisel \cite{Kneisel}, Lemma 2.2.1 and one also obtains a local representation for $D\Theta_h$. 
The latter yields the estimates by compactness of $\Sigma$.\hfill$\square$\\
\begin{Remark}\upshape \label{th_hanzawa_bem}
Let $h\in C^2(\Sigma)$ such that $\|h\|_\infty<a_0$. Then $\Gamma_h:=\Theta_h(\Sigma)$ is a connected and compact $C^2$-hypersurface as zero level set of
\[
\phi_h:\Sigma_a\rightarrow\R:x\mapsto d_\Sigma(x)-h(\Pi(x)),
\]
where $|\nabla\phi_h|\geq 1$ as one can show $\nabla\phi_h|_x\cdot\nu_\Sigma|_{\Pi(x)}=1$ for all $x\in\Sigma_a$ similar to Kneisel \cite{Kneisel}, p. 13.
Moreover, $\nu_{\Gamma_h}(x):=\nabla\phi_h(x)/|\nabla\phi_h(x)|$ is the outer unit normal to $\Gamma_h=\partial(\Theta_h(\Omega^+))$ at $x\in\Gamma_h$. Furthermore, because of $\phi_h(\Theta_h(x))=d_{\Sigma}(x)$ on $\overline{\Sigma_{a/3}}$, we get the identity $\nabla\phi_h|_{\Theta_h(x)}=D\Theta_h^{-\top}|_x\,\nu_\Sigma|_{\Pi(x)}$ for $x\in\overline{\Sigma_{a/3}}$ and this implies 
\begin{align}\label{eq_normale2}
D\Theta_h^{-\top}|_x\,\nu_\Sigma|_{\Pi(x)}\cdot\nu_\Sigma|_{\Pi(x)}=1\quad\text{ for all }x\in\overline{\Sigma_{a/3}}. 
\end{align} 
\end{Remark}

\subsection{Function Spaces for the Height Function and Modification $\tildeTheta_h$ of $\Theta_h$}\label{sec_fctraumh_modif_hanzawa}
In the following we specify, in which function spaces we consider height functions $h$. Let $p>n+2$ be fixed. We set
\[
X_0:=W^{1-\frac{1}{p}}_p(\Sigma)\quad\text{ and }\quad X_1:=W^{3-\frac{1}{p}}_p(\Sigma).
\]
Since $3-\frac{3}{p}-\frac{n-1}{p}>2$ we have
$
X_\gamma:=(X_0,X_1)_{1-\frac{1}{p},p}=W^{3-\frac{3}{p}}_p(\Sigma)\hookrightarrow C^2(\Sigma).
$ 
Moreover, we define $\E_1(T):=L^p(0,T;X_1)\cap W^1_p(0,T;X_0)$. Theorem \ref{th_einb_aus_trace_method} yields $\E_1(T)\hookrightarrow C^0([0,T];X_\gamma)$. Hence
\[
\Theta_h(x,t):=\Theta_{h(.,t)}(x)\quad \text{ for }(x,t)\in\R^n\times[0,T] 
\]
is well-defined for $h\in \E_1(T)\cap C^0([0,T];U)$ where $U:=\{h\in X_\gamma:\|h\|_\infty <a_0\}$. Let $R_0>0$ arbitrary and fixed from now on. It suffices to consider height functions in the space 
\[
V_T:=\{h\in\E_1(T)\cap C^0([0,T];U):\|h\|_{\E_1(T)}+\|h(0)\|_{X_\gamma}< R_0\}.
\]
The restriction due to the estimate is negligible as we are interested in local well-posedness, but it is needed to show suitable properties of a modification of the Hanzawa transformation $\Theta_h$ which we need for technical reasons: if we transform differential operators in space, the term $D\Theta_h^{-\top}$ always appears. In order to show appropriate regularity properties in space of $D\Theta_h^{-\top}$, it is complicated (but possible) to work with the Hanzawa transformation directly. However, it seems difficult to obtain suitable regularity in time for $D\Theta_h^{-\top}$ since derivatives of $h$ in space are present. It is much simpler to replace $h\circ\Pi$ in $\Theta_h$ by an extension $Eh$ of $h$, where for $\overline{a}:=\frac{3a}{4}$
\begin{align}\label{eq_fortsop}
E:X_0=W_p^{1-\frac{1}{p}}(\Sigma)\rightarrow W^1_p(\Sigma_{\overline{a}})\quad\text{ with } E\in\Lc(W_p^{k-\frac{1}{p}}(\Sigma),W^k_p(\Sigma_{\overline{a}})), k=1,2,3
\end{align}
is a suitable extension operator. Then we have $E\in\Lc(X_\gamma,\tilde{X}_\gamma)$, where we have set $\tilde{X}_\gamma:=(W^1_p(\Sigma_{\overline{a}}),W^3_p(\Sigma_{\overline{a}}))_{1-\frac{1}{p},p}$. Here because of $p>n+2$ it holds
\begin{align}\label{eq_fortsop2}
\tilde{X}_\gamma=W_p^{3-\frac{2}{p}}(\Sigma_{\overline{a}})\hookrightarrow C^2(\overline{\Sigma_{\overline{a}}}).
\end{align}

So we define the modification $\tildeTheta_h$ of the Hanzawa transformation $\Theta_h$ as follows: 
\begin{Definition}\label{th_def_hanzawa_mod}\upshape
Let $p>n+2,R_0>0$ as above and $T>0$. For $h_0\in U$ we set
\[
\tildeTheta_h(x):=x+\chi(\frac{d_\Sigma(x)}{a})Eh(x)\nu_\Sigma(\Pi(x))\quad\text{ for all }x\in\R^n
\]
and for $h\in V_T$ let $\tildeTheta_h(x,t):=\tildeTheta_{h(.,t)}(x)$ for all $(x,t)\in\R^n\times[0,T]$.
\end{Definition}

The following lemma shows that similar properties as for $\Theta_h$ also hold for the modification $\tildeTheta_h$ if the extension operator $E$ is chosen properly.
\begin{Lemma}\label{th_hanzawa_mod}
Let $\overline{a}=\frac{3a}{4}$ and $\varepsilon>0$. Then $E$ can be chosen such that for all $h\in X_\gamma$
\begin{align}\label{eq_fortsop_absch}
\|h\circ\Pi-Eh\|_{C^1(\overline{\Sigma_{\overline{a}}})}\leq \varepsilon\|h\|_{X_\gamma}.
\end{align}
Moreover, for $R>0$ there is an $\varepsilon=\varepsilon(R)>0$ such that with the extension operator $E$ with respect to $\varepsilon$ the following hold: $\tildeTheta_{h}:\Omega\rightarrow\Omega$ is for all $h\in U$ with $\|h\|_{X_\gamma}\leq R$ a $C^2$-diffeomorphism, the identity on $\R^n\textbackslash \Sigma_{2a/3}$ and we have
\[
\det D\tildeTheta_{h}\geq \tilde{c}>0\quad\text{ and }\quad D\tildeTheta_{h}^{-\top}|_x\nu_\Sigma|_{\Pi(x)}\cdot \nu_\Sigma|_{\Pi(x)}\geq \frac{1}{2}\quad \text{ for all }x\in\overline{\Sigma_{a/3}}.
\]
\end{Lemma}
\begin{Remark}\label{th_hanzawa_mod_bem}\upshape
Theorem \ref{th_einb_aus_trace_method} implies $\E_1(T)\hookrightarrow C^0([0,T];X_\gamma)$ and the embedding constant is bounded independent of $T>0$ if we add the term $\|h(0)\|_{X_\gamma}$ in the $\E_1(T)$-norm. Now we can do the following: we choose $R=\overline{R_0}>R_0$ in Lemma \ref{th_hanzawa_mod} (and so $E$) such that
\[
h(t)\in U_0:=\{h_0\in U:\|h_0\|_{X_\gamma}<\overline{R_0}\}
\] 
holds for all $h\in V_T$ and $0\leq t\leq T$. In particular for $h\in V_T$ the last part in Lemma \ref{th_hanzawa_mod} is valid for $\tildeTheta_h(.,t)$ instead of $\tildeTheta_h$. Here $0\leq t\leq T$ and $0<T<\infty$ are arbitrary. 
\end{Remark}
\textit{Proof of Lemma \ref{th_hanzawa_mod}.} Similar to A. and Wilke \cite{Abels_Mullins_Sekerka} one can construct a new operator $E$ that fulfils \eqref{eq_fortsop_absch} by scaling $E$ in normal direction of $\Sigma$. Now let $R>0$ arbitrary. First of all, $\tildeTheta_{h}$ is well-defined for all $h\in X_\gamma$, a $C^2$-mapping and the identity on $\R^n\textbackslash\Sigma_{2a/3}$, as the second part is cut off properly. Since $\tildeTheta_h=\Theta_h$ on $\Sigma$, we can use the properties of $\Theta_{h}$. For all $h\in U$ with $\|h\|_{X_\gamma}\leq R$ it holds $\|h\circ\Pi-Eh\|_{C^1(\overline{\Sigma_{\overline{a}}})}\leq \varepsilon R$ due to \eqref{eq_fortsop_absch} and therefore
\begin{align}\label{eq_bew_hanzawa_differenz}
\|D\Theta_{h}-D\tildeTheta_{h}\|_\infty\leq\|\chi(\frac{d_\Sigma(.)}{a})\nu_\Sigma(\Pi(.))(h\circ\Pi-Eh)\|_{C^1(\overline{\Sigma_{\overline{a}}})}\leq C\varepsilon R.
\end{align}
Since $\{D\Theta_{h}(x):x\in\R^n,h\in U, \|h\|_{X_\gamma}\leq R\}$ is contained in a compact set $K\subseteq GL(n,\R)$, we obtain for $\varepsilon(R)>0$ small that $\det D\tildeTheta_{h}\geq \frac{c}{2}>0$. In particular $D\tildeTheta_{h}$ is invertible and the inverse mapping theorem implies that $\tilde{\Theta}_h:\Omega\rightarrow\Omega$ is a $C^2$-diffeomorphism if we show bijectivity of $\tildeTheta_h:\R^n\rightarrow\R^n$. This can be done for a suitable choice of $\varepsilon(R)$ as in A. and Wilke \cite{Abels_Mullins_Sekerka}, p. 49 using Lemma \ref{th_hanzawa} and \eqref{eq_bew_hanzawa_differenz}.

It remains to show the second estimate. Lemma \ref{th_hanzawa} and \eqref{eq_bew_hanzawa_differenz} yield because of $D\tildeTheta_h=D\Theta_h(I-D\Theta_h^{-1}(D\Theta_h-D\tildeTheta_h))$ together with a Neumann series argument
\[
|D\tildeTheta_{h}^{-\top}(x)-D\Theta_{h}^{-\top}(x)|\leq \tilde{C}(R)\varepsilon
\]
for all $x\in\R^n$ and $h\in U$ with $\|h\|_{X_\gamma}\leq R$ if $\varepsilon=\varepsilon(R)>0$ is sufficiently small. Equation \eqref{eq_normale2} in Remark \ref{th_hanzawa_bem} implies
\[
|D\tildeTheta_{h}^{-\top}|_x\nu_\Sigma|_{\Pi(x)}\cdot\nu_\Sigma|_{\Pi(x)}-1|=|(D\tildeTheta_{h}^{-\top}-D\Theta_{h}^{-\top})|_x\nu_\Sigma|_{\Pi(x)}\cdot\nu_\Sigma|_{\Pi(x)}|\leq \tilde{C}(R)\varepsilon \leq \frac{1}{2}
\]
for all $x\in\overline{\Sigma_{a/3}}$ if $\varepsilon(R)>0$ is chosen suitably.\hfill$\square$\\

\subsection{The Transformed Equations}\label{sec_transf_glg}
With $\tildeTheta_h$ we transform the Navier-Stokes/mean curvature flow system \eqref{eq_nsmc_1}-\eqref{eq_nsmc_8} to $\Sigma$ and $\Omega^\pm$, respectively. Therefore we fix the height function $h\in V_T$ and let $v:\Omega\times(0,T)\rightarrow\R^n$ be a vector-valued and $\tilde{p}:\Omega\times(0,T)\rightarrow\R$ be a scalar-valued function. Then we define $u(x,t):=v(\tildeTheta_h(x,t),t)$ and $\tilde{q}(x,t):=\tilde{p}(\tildeTheta_h(x,t),t)$ for $(x,t)\in\Omega\times(0,T)$. Moreover, let
\[
A(h):=D\tildeTheta_h^{-\top},\quad\nabla_h:=A(h)\nabla,\quad \nabla_hu:=(\nabla_hu_k^\top)_{k=1}^n\:\text{ and }\:\diverg_hu:=\Tr\nabla_hu.
\]
Let $t\in[0,T]$. We set $\Gamma_h(t):=\Gamma_{h(.,t)}$ and denote the outer unit normal by $\nu_{\Gamma_h(t)}$ and by $H_{\Gamma_h(t)}$ the corresponding mean curvature. Furthermore, we define
\[
K(h):\Sigma\times[0,T]\rightarrow\R:(x,t)\mapsto H_{\Gamma_h(t)}|_{\Theta_h(x,t)}.
\]
Applying the chain rule to $\phi_{h(.,t)}\circ\tildeTheta_h(.,t)$ with $\phi_{h(.,t)}$ as in Remark \ref{th_hanzawa_bem} implies 
\[
\nu_{h}|_{(x,t)}:=\nu_{\Gamma_h(t)}|_{\tildeTheta_h(x,t)}=\frac{A(h)|_{(x,t)}\nu_{\Sigma}|_x}{|A(h)|_{(x,t)}\nu_{\Sigma}|_x|}\quad\text{ for all }(x,t)\in\Sigma\times(0,T).
\] 
Moreover, we have
$V_{\Gamma_h(t)}|_{\tildeTheta_h(x,t)}=\partial_t(\Theta_h(x,t))\cdot\nu_{\Gamma_h(t)}|_{\Theta_h(x,t)}=\partial_th(x,t)\,\nu_\Sigma|_x\cdot\nu_h|_{(x,t)}$ for the normal velocity. Since $\tildeTheta_h$ is the identity in a neighbourhood of $\partial\Omega$, the no-slip condition on $v$ is retained. Therefore \eqref{eq_nsmc_1}-\eqref{eq_nsmc_8} becomes
\begin{alignat}{2}
\partial_tu+u\cdot\nabla_hu-\mu^\pm\diverg_h\nabla_h u+\nabla_h \tilde{q}&=\nabla_hu\,\partial_t\tildeTheta_h &\qquad&\text{ in }\Omega^\pm\times(0,T),\label{eq_transf_glg1}\\
\diverg_hu&=0 &\qquad&\text{ in }\Omega^\pm\times(0,T),\label{eq_transf_glg2}\\
-\lsprung  T_h(u,\tilde{q})\rsprung\,\nu_h&=\sigma K(h)\,\nu_h &\qquad&\text{ on }\Sigma\times(0,T),\label{eq_transf_glg3}\\
\lsprung u\rsprung &=0 &\qquad&\text{ on }\Sigma\times(0,T),\label{eq_transf_glg4}\\
u|_{\partial\Omega}&=0 &\qquad&\text{ on }\partial\Omega\times(0,T),\label{eq_transf_glg5}\\
u|_{t=0}&=v_0\circ \tildeTheta_{h_0} &\qquad&\text{ in }\Omega,\label{eq_transf_glg6}\\
\partial_t h-a_1(h)\cdot u|_\Sigma &=a_2(h)K(h) &\qquad&\text{ on }\Sigma\times(0,T),\label{eq_transf_glg7}\\
h|_{t=0}&=h_0 &\qquad&\text{ on }\Sigma,\label{eq_transf_glg8}
\end{alignat}
where we set $T_h(u,\tilde{q}):=2\mu^\pm\text{sym}(\nabla_hu)-\tilde{q}I$ as well as
\[
a_1(h):=\frac{A(h)\nu_\Sigma}{A(h)\nu_\Sigma\cdot\nu_\Sigma}\quad\text{ and }\quad a_2(h):=\frac{|A(h)\nu_\Sigma|}{A(h)\nu_\Sigma\cdot\nu_\Sigma}.
\]
Moreover $u\cdot\nabla_hu\equiv (u\cdot\nabla_h)u$ and $\diverg_h$ is applied row-wise to $\nabla_hu$ in \eqref{eq_transf_glg1}. For this system we show local well-posedness. As a preparation, in the next two Subsections \ref{sec_A_h} and \ref{sec_mittl_krg} we show properties of $\tildeTheta_h, A(h)=D\tildeTheta_h^{-\top}$ and of the transformed mean curvature $K(h)$, respectively.

\subsection{Properties of $\tildeTheta_h$ and $D\tildeTheta_h^{-\top}$}\label{sec_A_h} 
For later use we must know the regularity of $\tildeTheta_h$ and $A(h)=D\tildeTheta_h^{-\top}$. This is provided by
\begin{Lemma}\label{th_A}
Let $U_0$ as in Remark \ref{th_hanzawa_mod_bem} and $0<T\leq T_0$. Then $\tildeTheta_.\in BC^1(U_0;C^2(\overline{\Omega})^n)$, $A\in BC^1(U_0;C^1(\overline{\Omega})^{n\times n})$ and $\tildeTheta_.\in BC^1(V_T;\tilde{X}_T), A\in BC^1(V_T;X_T)$, where
\begin{align*}
\tilde{X}_T&:=C^\tau([0,T];C^2(\overline{\Omega})^{n})\times C^{\frac{1}{2}+\tau}([0,T];C^1(\overline{\Omega})^{n})\cap W^1_p(0,T;W^1_p(\Omega)^{n}),\\
X_T&:=C^\tau([0,T];C^1(\overline{\Omega})^{n\times n})\times C^{\frac{1}{2}+\tau}([0,T];C^0(\overline{\Omega})^{n\times n})\cap W^1_p(0,T;L^p(\Omega)^{n\times n})
\end{align*}
for some $\tau>0$. In the time-dependent case the mappings and their derivatives are bounded by a $C>0$ independent of $0<T\leq T_0$ if we add $\|h(0)\|_{X_\gamma}$ in the $\E_1(T)$-norm.
\end{Lemma}
\textit{Proof.} Let $h\in U_0$ or $h\in V_T$. First, we show the properties of $\tildeTheta_h$. It holds
\begin{align}\label{eq_bew_Ah_phi}
\tildeTheta_h=\textup{id}+\phi Eh\quad\text{ with }\phi:=\chi(d_\Sigma(.)/a)\,\nu_\Sigma\circ\Pi\in C^\infty_0(\Sigma_{3a/4})^n.
\end{align}
Here $E$ is the extension operator from Remark \ref{th_hanzawa_mod_bem} that induces $E\in\Lc(X_\gamma,\tilde{X}_\gamma)$ with $\tilde{X}_\gamma\hookrightarrow C^2(\overline{\Sigma_{\overline{a}}})$ by \eqref{eq_fortsop2} and $E\in\Lc(\E_1(T),\tilde{\E}_1(T))$ with
\[ \tilde{\E}_1(T):=L^p(0,T;W^3_p(\Sigma_{\overline{a}}))\cap W^1_p(0,T;W^1_p(\Sigma_{\overline{a}})),
\]
where the operator norm is bounded by a $C>0$ independent of $T>0$. We need suitable embeddings for $\tilde{\E}_1(T)$: Theorem \ref{th_einb_E_1T_in_hölder} yields for $0<\theta<1-\frac{1}{p}$
\begin{align}\label{eq_bew_Ah_0}
\tilde{\E}_1(T)\hookrightarrow C^{1-\theta-\frac{1}{p}}([0,T];B^{1+2\theta}_{p,p}(\Sigma_{\overline{a}}))
\end{align}
and the embedding constant is bounded independent of $T>0$ if we add the term $\|\tilde{h}(0)\|_{\tilde{X}_\gamma}$ in the $\tilde{\E}_1(T)$-norm. For $\theta=\frac{n}{2p}+\frac{k-1}{2}+\delta, k=1,2$ and $\delta>0$ small we obtain because of $p>n+2$
\[
1+2\theta-\frac{n}{p}=1+k-1+2\delta>k\quad\text{ and }\quad 1-\frac{1}{p}-\theta=1-\frac{n+2}{2p}-\frac{k-1}{2}-\delta>1-\frac{k}{2}.
\]
Embedding theorems and \eqref{eq_bew_Ah_0} imply $\tilde{\E}_1(T)\hookrightarrow C^\tau([0,T];C^2(\overline{\Sigma_{\overline{a}}}))\cap C^{\frac{1}{2}+\tau}([0,T];C^1(\overline{\Sigma_{\overline{a}}}))$ for some $\tau>0$. Since $\tildeTheta_h$ is affine linear in $h$, we get the properties of $\tildeTheta_h$.

Now we deduce the assertions for $A(h)=D\tildeTheta_h^{-\top}$. The first part yields analoguous properties for $D\tildeTheta_h$. Moreover, the space $X_T$ is an algebra with pointwise multiplication and a product estimate holds because this is valid for the Hölder spaces and we also have Lemma \ref{th_W_C0}. Additionally, Lemma \ref{th_hanzawa_mod} implies $\det D\tildeTheta_h\geq c>0$ for all $h\in U_0$ and $h\in V_T$, respectively. The inverse formula, Lemma \ref{th_hölder_ableiten} and Lemma \ref{th_W_C0} give
\begin{align}\label{eq_bew_Ah_1}
A(h)\in C^1(\overline{\Omega})^{n\times n}\text{ for all }h\in U_0\quad\text{ and }\quad A(h)\in X_T\text{ for all }h\in V_T.
\end{align}
Furthermore, the mappings are bounded and in the time-dependent case bounded by a $C>0$ independent of $0<T\leq T_0$ if we add $\|h(0)\|_{X_\gamma}$ in the $\E_1(T)$-norm. The implicit function theorem applied to
\[
G:Z\times X\rightarrow X: (h,B)\mapsto D\tildeTheta_h^\top B-I
\]
at $(h,A(h))$ for $(Z,X)=(U_0,C^1(\overline{\Omega})^{n\times n})$ or $(V_T,X_T)$ implies $A\in C^1(U_0;C^1(\overline{\Omega})^{n\times n})$ and $A\in C^1(V_T;X_T)$, respectively, with
\[
\frac{d}{dh}A(h)(.)=-A(h)\frac{d}{dh}(D\tildeTheta_h^{\top})(.)A(h).
\]
From the properties of $D\tildeTheta_h$ and \eqref{eq_bew_Ah_1} the claim follows.\hfill$\square$

As a consequence we obtain
\begin{Corollary}\label{th_A_kor}
Let $U_0$ be as in Remark \ref{th_hanzawa_mod_bem}, $0<T\leq T_0$ and $\tau>0$ as in Lemma \ref{th_A}. Then $\nu_h, a_1(h)$ and $a_2(h)$ as in Subsection \ref{sec_transf_glg} are well-defined for $h\in U_0$ and $h\in V_T$, respectively, and have the regularity
\[
BC^1(U_0;C^1(\Sigma))\quad\text{ and }\quad BC^1(V_T;C^\tau([0,T];C^1(\Sigma))\cap C^{\frac{1}{2}+\tau}([0,T];C^0(\Sigma))).
\]
In the time-dependent case the mappings and their derivatives are bounded by a $C>0$ independent of $0<T\leq T_0$ if we add the term $\|h(0)\|_{X_\gamma}$ in the $\E_1(T)$-norm.
\end{Corollary}
\textit{Proof.} Lemma \ref{th_hanzawa_mod} and Remark \ref{th_hanzawa_mod_bem} yield $|A(h)\nu_\Sigma|_\Pi|\geq \frac{1}{2}$ on $\overline{\Sigma_{a/3}}$. Thus we can extend $\nu_\Sigma$ by $\nu_\Sigma|_\Pi$ smoothly to $\overline{\Sigma_{a/3}}$ and obtain the assertion by Lemma \ref{th_A}.\hfill$\square$\\

\subsection{Mean Curvature}\label{sec_mittl_krg}
For $h\in U$ we denote by $H_{\Gamma_h}$ the mean curvature of $\Gamma_h=\Theta_h(\Sigma)$ with respect to $\nu_{\Gamma_h}$. This subsection is devoted to the transformed mean curvature $K(h):=H_{\Gamma_h}\circ\Theta_h:\Sigma\rightarrow\R$. In particular we show that $K(h)$ has a quasilinear structure in terms of $h$ and the principal part for $h=0$ is given by the Laplace-Beltrami operator $\Delta_\Sigma$.
\begin{Lemma}\label{th_mittl_krg}
Let $p>n+2$. Then there are $P\in C^1(U;\Lc(X_1,X_0))$ and $Q\in C^1(U;X_0)$ with $P(0)=\Delta_\Sigma$ and $K(h)=P(h)h+Q(h)$ for all $h\in U\cap X_1$.
\end{Lemma}
\textit{Proof.} Let $\varphi_l:V_l\rightarrow U_l\subseteq\Sigma$ for $l=1,...,N$ be suitable parametrizations of $\Sigma$ with $\Sigma=\bigcup_{l=1}^N U_l$ and $V_l\subseteq\R^{n-1}$. Then $\tilde{\varphi}_l:=(\varphi_l,\textup{id}):V_l\times(-a,a)\rightarrow U_l\times(-a,a)$
are appropriate parametrizations of $\Sigma\times(-a,a)$. Now fix $l$. The Euclidean metric $g^\textup{eucl}$ on $\Sigma_a$ induces via the diffeomorphism $X$ from Subsection \ref{sec_hanzawa} a Riemannian metric $g^X$ on $\Sigma\times(-a,a)$. We denote by $w_{ij}, \Gamma^k_{ij}$ for $i,j,k=1,...,n$ the local representation of $g^X$ and the Christoffel symbols with respect to $\tilde{\varphi}_l$, respectively, and $w^{ij}$ is as usual. Moreover, we set $w_{ij}(h)|_s:=w_{ij}|_{(s,h(s))}$ for $s\in U_l$. $\Gamma_{ij}^k(h)$ and $w^{ij}(h)$ are defined analogously. As in Escher and Simonett \cite{EscherSimonett_MullinsSekerka}, proof of Lemma 3.1 on p. 274 ff. and Remark 3.2 on p. 277 it follows that
\[
K(h)|_{U_l}=P_l(h)h+Q_l(h)
\]
where $P_l$ and $Q_l$ are altered by a factor $-(n-1)$ here. Additionally let $p_{jk}, p_j, q, l_.$ be as in Escher and Simonett \cite{EscherSimonett_MullinsSekerka} modified by the same factor. Then $P_l(0)$ is a local representation of $\Delta_\Sigma$. One can show that the $w_{ij}|_{(s,r)}$ and $\partial_i\partial_jX\cdot\partial_mX|_{(s,r)}$ for $i,j,m=1,...,n$ are at most quadratic in $r\in(-a,a)$ with smooth coefficients in $s\in \overline{U_l}$. Because of $\|h\|_\infty<a_0<\frac{a}{4}$ and $w_{ij}(h)=w_{ij}(.,h(.))$ for $i,j=1,...,n$, a compactness argument yields
\begin{align}\label{eq_det_w_ij_absch}
\text{det}((w_{ij}(h))_{i,j=1}^n)\geq c>0
\end{align} 
for a $c>0$ independent of $h$.

Using $X_\gamma\hookrightarrow C^2(\Sigma)$, the product rule for the Fréchet-derivative, the inverse formula as well as
Lemma \ref{th_hölder_ableiten} for the local representations in Escher and Simonett \cite{EscherSimonett_MullinsSekerka} implies
\[
p_{jk}(.)\circ\varphi_l,\,p_i(.)\circ\varphi_l,\, q(.)\circ\varphi_l\in C^1(U;C^1(\overline{V_l}))\quad\text{ for }i,j,k=1,...,n-1.
\]
Now the claim follows utilizing a suitable partition of unity.\hfill$\square$

Furthermore, we need properties of $K$, when we insert a time-dependent height function $h\in V_T$. This gives us
\begin{Lemma}\label{th_mittl_krg2}
Let $p>n+2, 1<q<p$ and $0<T\leq T_0$. Then it holds
\[
K\in BC^1(V_T;W^{1-\frac{1}{q},\frac{1}{2}(1-\frac{1}{q})}_{q}(\Sigma_T)).
\]
Moreover, $K$ and its derivative are bounded by a constant $C(p,q,T_0)>0$ when we add the term $\|h(0)\|_{X_\gamma}$ in the $\E_1(T)$-norm.
\end{Lemma}
\textit{Proof.} We use the same notation as in the proof of Lemma \ref{th_mittl_krg}. Then 
\begin{align}\label{eq_bew_mc0}
K(h)|_{U_l}=\sum_{\alpha\in\N_0^{n-1},|\alpha|\leq2}a_\alpha(.,h,\nabla_sh)\partial_s^\alpha h,
\end{align}
where the $a_\alpha$ correspond to $p_{ij},p_j$ and $q$, respectively, and $\partial_s^\alpha$ is the $\alpha$-th derivative with respect to the canonical coordinates. 

First we estimate the $a_\alpha$ in suitable Hölder spaces. As in the proof of Lemma \ref{th_A} one can show $\E_1(T)\hookrightarrow C^{\tau}([0,T];C^2(\Sigma))\cap C^{\frac{1}{2}+\tau}([0,T];C^1(\Sigma))$ for $\tau>0$ small because of $p>n+2$ and the embedding constant is bounded independent of $T>0$ if we add the term $\|h(0)\|_{X_\gamma}$ in the $\E_1(T)$-norm. Since the $w_{ij}$ for $i,j=1,...,n$ are at most quadratic in $h$ with smooth coefficients in $\overline{U_l}$ product estimates in Hölder spaces yield
\[
w_{ij}(.)\circ\varphi_l\in BC^1(V_T;C^\tau([0,T];C^2(\overline{V_l}))\cap C^{\frac{1}{2}+\tau}([0,T];C^1(\overline{V_l})))\quad\text{ for }i,j=1,...,n.
\]
Now \eqref{eq_det_w_ij_absch} implies\footnote{ In fact, the same regularity as for the $w_{ij}(.)\circ\varphi_l$ should be obtained. But this is not needed.} using the inverse formula, product estimates and Lemma \ref{th_hölder_ableiten} that
\[
w^{ij}(.)\circ\varphi_l\in BC^1(V_T;C^\tau([0,T];C^1(\overline{V_l}))\cap C^{\frac{1}{2}+\tau}([0,T];C^0(\overline{V_l})))\quad\text{ for }i,j=1,...,n.
\] 
Hence this also holds for $\Gamma_{ij}^k(.)\circ\varphi_l$ and $(l_.)^r$, where $i,j,k=1,...,n$ and $r\in\R$. Using the mapping properties in the local representations, we obtain for $\alpha\in\N_0^{n-1}, |\alpha|\leq 2$
\begin{align}\label{eq_bew_mc_a_alpha}
a_\alpha\in BC^1(V_T;C^\tau([0,T];C^1(\overline{V_l}))\cap C^{\frac{1}{2}+\tau}([0,T];C^0(\overline{V_l}))).
\end{align}
Furthermore, all above terms and their derivatives are bounded by a $C>0$ independent of $T>0$ if we add $\|h(0)\|_{X_\gamma}$ in the $\E_1(T)$-norm. So this also holds for the $a_\alpha$. The exponent $q$ does not appear here.

To show the desired properties of $K$ we use \eqref{eq_bew_mc0}. Therefore we need another embedding for $\E_1(T)$ to estimate $\partial_s^\alpha h$ suitably. Theorem \ref{th_intungl_einb_L_B} and a well-known interpolation inequality for real interpolation spaces (see Corollary 1.7 in Lunardi \cite{Lunardi_Interpolation}) yield
\[
\E_1(T)\hookrightarrow B^{1-\theta}_{p,\infty}(0,T;(X_0,X_1)_{\theta,p})
\] 
for $\theta\in(0,1)$ and the embedding constant is bounded independent of $T>0$. We set $\theta=\frac{1}{2}(1+\frac{1}{p})+\varepsilon$ for $\varepsilon>0$ small. Then for $\varepsilon=\varepsilon(p,q)>0$ small it holds that
\[
1-\frac{1}{p}+2\theta=2+2\varepsilon\quad\text{ and }\quad 1-\theta=\frac{1}{2}(1-\frac{1}{p})-\varepsilon>\frac{1}{2}(1-\frac{1}{q}).
\]
We infer $
\E_1(T)\hookrightarrow B_{p,\infty}^{\frac{1}{2}(1-\frac{1}{p})-\varepsilon}(0,T;W^2_p(\Sigma))$ and, again, the embedding constant is bounded independent of $T>0$. This implies
\[
\partial_{x_i}(.\circ\varphi_l),\, \partial_{x_j}\partial_{x_i}(.\circ\varphi_l)\in BC^1(V_T;B_{p,\infty}^{\frac{1}{2}(1-\frac{1}{p})-\varepsilon}(0,T;L^p(V_l))\cap L^p(0,T;W_p^{1-\frac{1}{p}}(V_l)))
\]
for $i,j=1,...,n-1$ and the mappings and their derivatives are bounded independent of $T>0$. Since pointwise multiplication is a product on $C^1(\overline{V_l})\times W^{1-\frac{1}{p}}_p(V_l)\rightarrow W^{1-\frac{1}{p}}_p(V_l)$ in sense of Definition \ref{th_def_produkt}, using Lemma \ref{th_hölder_prod_sob_bes} and \eqref{eq_bew_mc0}-\eqref{eq_bew_mc_a_alpha} yields
\[
K(.)\circ\varphi_l \in BC^1(V_T;B_{p,\infty}^{\frac{1}{2}(1-\frac{1}{p})-\varepsilon}(0,T;L^p(V_l))\cap L^p(0,T;W_p^{1-\frac{1}{p}}(V_l)))
\]
where $K(.)\circ\varphi_l$ and the derivative are bounded by a constant $C>0$ independent of $T>0$. Now we use the embeddings $L^p(V_l)\hookrightarrow L^q(V_l)$ and $W_p^{1-\frac{1}{p}}(V_l)\hookrightarrow W_q^{1-\frac{1}{q}}(V_l)$ as well as 
\[
B_{p,\infty}^{\frac{1}{2}(1-\frac{1}{p})-\varepsilon}(0,T;L^p(V_l))\hookrightarrow B_{q,\infty}^{\frac{1}{2}(1-\frac{1}{p})-\varepsilon}(0,T;L^p(V_l))
\] 
where the latter can be shown by using Hölder's inequality and the embedding constant is bounded by $C(p,q,T_0)>0$. We infer from Simon \cite{Simon}, Corollary 15 that
\[
B_{q,\infty}^{\frac{1}{2}(1-\frac{1}{p})-\varepsilon}(0,T;L^q(V_l))\hookrightarrow W_q^{\frac{1}{2}(1-\frac{1}{q})}(0,T;L^q(V_l))
\]
and the embedding constant is also bounded by a $C(p,q,T_0)>0$. The claim follows by using a suitable partition of unity.\hfill$\square$\\

\section{Two-Phase Navier-Stokes System for Given Time-Dependent Interface}\label{sec_NavSt_teil}
We consider the Navier-Stokes part \eqref{eq_transf_glg1}-\eqref{eq_transf_glg6} for a fixed time-dependent height function $h\in V_T$ (or equivalently for a given family of hypersurfaces $\{\Gamma_{h}(t)\}_{t\in(0,T)}$) and show unique solvability and some properties of the solution operator. Let as before $p>n+2, R_0>0$ be fixed as well as $X_0,X_1,X_\gamma,U,\E_1(T)$ and $V_T$ for $T>0$ as in Subsection \ref{sec_fctraumh_modif_hanzawa}. Moreover, let $2<q<3$ and $h\in V_T$ with $h_0:=h(0)$. Then we set
\[
Y_{\gamma,h_0}:= \left\{ v\in L^p(\Omega):v|_{\tildeTheta_{h_0}(\Omega^\pm)}\in W^{2-\frac{2}{q}}_{q}(\tildeTheta_{h_0}(\Omega^\pm))\right\}^n\cap W^1_{q,0}(\Omega)^n\quad\text{ and }\quad Y_\gamma:=Y_{\gamma,0}.
\]
Now let $v_0\in Y_{\gamma,h_0}\cap L^q_{\sigma}(\Omega)$. In order to write the Navier-Stokes part \eqref{eq_transf_glg1}-\eqref{eq_transf_glg6} as an abstract fixed-point equation, we rearrange the equations such that the left hand side is the same as for the linear Stokes system in Appendix \ref{sec_lin_stokes}. This yields
\begin{alignat}{2}
\partial_tu-\mu^\pm\Delta u+\nabla \tilde{q}&=a^\pm(h;D_x)(u,\tilde{q})-u\cdot\nabla_hu+\nabla_hu\,\partial_t\tildeTheta_h &\quad&\text{ in }\Omega^\pm\times(0,T),\label{eq_navstk_transf1}\\
\diverg\,u&= \Tr(\nabla u-\nabla_h u)=:g(h)u &\quad&\text{ in }\Omega^\pm\times(0,T),\label{eq_navstk_transf2}\\
\lsprung T(u,\tilde{q})\rsprung\,\nu_\Sigma&=t^\pm(h;D_x)(u,\tilde{q})-\sigma K(h)\nu_h &\quad&\text{ on }\Sigma\times(0,T),\label{eq_navstk_transf3}\\
\lsprung u\rsprung&=0 &\quad&\text{ on }\Sigma\times(0,T),\label{eq_navstk_transf4}\\
u|_{\partial\Omega}&=0 &\quad&\text{ on }\partial\Omega\times(0,T),\label{eq_navstk_transf5}\\
u|_{t=0}&=v_0\circ\tildeTheta_{h_0} &\quad&\text{ in }\Omega,\label{eq_navstk_transf6}
\end{alignat}
where $A(h)=D\tildeTheta_h^{-\top}, \nabla_h=A(h)\nabla$ and $K(h),\nu_h$ are as in Subsection \ref{sec_transf_glg} as well as
\begin{align*}\begin{split}
a^\pm(h;D_x)(u,\tilde{q})&:=\mu^\pm(\diverg_h\nabla_hu-\Delta u)+(\nabla-\nabla_h)\tilde{q},\\
t^\pm(h;D_x)(u,\tilde{q})&:=\lsprung 2\mu^\pm \text{sym}(\nabla u)-\tilde{q}I\rsprung\,(\nu_\Sigma-\nu_h)+\lsprung 2\mu^\pm\text{sym}(\nabla u-\nabla_hu)\rsprung\,\nu_h.\end{split}
\end{align*}
We introduce the space $Y_{T,h_0,v_0}:=\{(u,\tilde{q})\in Y_T :u|_{t=0}=v_0\circ\tildeTheta_{h_0}\}$ for the transformed velocity and pressure $(u,\tilde{q})$, where $Y_T:=Y_T^1\times Y_T^2$ and with $\Omega_0:=\Omega^+\cup\Omega^-$ we define
\begin{align}\begin{split}\label{eq_Y_T}
Y_T^1&:= W^1_q(0,T;L^q(\Omega))^n\cap L^q(0,T;W^2_q(\Omega_0)\cap W^1_{q,0}(\Omega))^n,\\
Y_T^2&:=\{\tilde{q}\in L^q(0,T;W^1_{q,(0)}(\Omega_0)):\lsprung \tilde{q} \rsprung\in 
W_{q}^{1-\frac{1}{q},\frac{1}{2}(1-\frac{1}{q})}(\Sigma_T)\}.\end{split}
\end{align}

Then the existence result for the Navier-Stokes part reads as follows:
\begin{Theorem}\label{th_navstk_exsatz}
Let $p>n+2$ and $2<q<3$ with $1+\frac{n+2}{p}>\frac{n+2}{q}$. Then for $R>0$ there are $0<\varepsilon=\varepsilon(R)<R_0$ and $T_0=T_0(R)>0$ such that for all $0<T\leq T_0$ and 
\[
h\in V_{T,\varepsilon}:=\{ h\in V_T:  \|h(0)\|_{X_\gamma}\leq\varepsilon \}, \quad v_0\in Y_{\gamma,h(0)}\cap L^q_{\sigma}(\Omega)\text{ with } \|v_0\|_{Y_{\gamma,h(0)}}\leq R
\]
the Navier-Stokes part \eqref{eq_navstk_transf1}-\eqref{eq_navstk_transf6} has a unique solution $F_T(h,v_0):=(u,\tilde{q})(h)\in Y_{T,h(0),v_0}$ and it holds that $\|F_T(h,v_0)\|_{Y_T}\leq C_R$ as well as $F_T(h,v_0)|_{[0,\tilde{T}]}=F_{\tilde{T}}(h|_{[0,\tilde{T}]},v_0)$ for all $0<\tilde{T}\leq T$.

Additionally, for $h_j\in V_{T,\varepsilon}$ with $h_0^j:=h_j(0)$ as well as $v_0^j\in Y_{\gamma,h_j(0)}\cap L^q_{\sigma}(\Omega)$ with $\|v_0^j\|_{Y_{\gamma,h_j(0)}}\leq R$ and $\tilde{v}_0^j:=v_0^j\circ\tildeTheta_{h_0^j}$ for $j=1,2$ we have the estimate
\[
\|F_T(h_1,v_0^1)-F_T(h_2,v_0^2)\|_{Y_T}\leq \tilde{C}_R(\|h_1-h_2\|_{\E_1(T)}+\|h_0^1-h_0^2\|_{X_\gamma}+\|\tilde{v}_0^1-\tilde{v}_0^2\|_{Y_\gamma}).
\]
\end{Theorem}

For the proof we reformulate \eqref{eq_navstk_transf1}-\eqref{eq_navstk_transf6} as an abstract fixed-point equation
\begin{align}\label{eq_fpgl_stokes}
w=L^{-1}G(w;h,v_0)\quad\text{ for }w\in Y_{T,h(0),v_0}, 
\end{align} 
where for $w=(u,\tilde{q})\in Y_T$ we define
\begin{align*}
Lw:=
\begin{pmatrix}
\partial_tu-\mu^\pm\Delta u+\nabla\tilde{q}\\
\diverg\,u\\
\lsprung T(u,\tilde{q})\rsprung\,\nu_\Sigma\\
u|_{t=0}
\end{pmatrix},\,\,
G(w;h,v_0):=
\begin{pmatrix}
a^\pm(h;D_x)w-u\cdot\nabla_hu+\nabla_hu\,\partial_t\tildeTheta_h\\
g(h)u-\frac{1}{|\Omega|}\int_\Omega g(h)u\,dx\\
t^\pm(h;D_x)w-\sigma K(h)\nu_h\\
v_0\circ\tildeTheta_{h(0)}
\end{pmatrix}.
\end{align*}
Lemma \ref{th_stokes_lin_notw} and Theorem \ref{th_stokes_lin} imply that $L:Y_T\rightarrow Z_T$ is an isomorphism, where 
\[
Z_T:=\{(f,g,a,u_0)\in\tilde{Z}_T\,:\,\diverg\, u_0=g|_{t=0}\text{ in }W^{-1}_{q,(0)}(\Omega)\}
\]
with $\tilde{Z}_T:=\tilde{Z}_T^1\times \tilde{Z}_T^2\times (\tilde{Z}_T^3)^n\times Y_\gamma$ and the spaces
\begin{align*}
\tilde{Z}_T^1&:=L^q(0,T;L^q(\Omega))^n,\quad \tilde{Z}_T^3:= W_{q}^{1-\frac{1}{q},\frac{1}{2}(1-\frac{1}{q})}(\Sigma_T),\\
\tilde{Z}_T^2&:=\{ g\in L^q(0,T;W^1_{q,(0)}(\Omega_0))\cap W^1_q(0,T;W^{-1}_{q,(0)}(\Omega)):
\tr_\Sigma (g|_{\Omega^+})\in \tilde{Z}_T^3\}.
\end{align*} 
Here $Y_T,Z_T$ and $\tilde{Z}_T$ are Banach spaces with canonical norms, see Remark \ref{th_bem_stokes_lin}.

At this point we see why we have to apply\footnote{ A. and Wilke \cite{Abels_Mullins_Sekerka} transform with $\tildeTheta_h\circ\tildeTheta_{h_0}^{-1}$. Then for fixed $h_0$ one could omit the dependence on $v_0$. But here we have chosen the other strategy for the following reasons: On one hand it is difficult to obtain suitable properties of $\tildeTheta_{h_0}^{-1}$ and on the other hand $\Gamma_0:=\tildeTheta_{h_0}(\Sigma)$ would have the role of $\Sigma$ in Theorem \ref{th_stokes_lin} but one needs a priori $C^{2,1}$-regularity for the interface there, see Shimizu \cite{Shimizu}.} the fixed-point argument in the space $Y_{T,h_0,v_0}$ which depends on $h_0$ and $v_0$. Without the additional condition for $u|_{t=0}$ in $Y_{T,h(0),v_0}$ the compatibility condition for $G(w;h,v_0)$ in $Z_T$ would not be fulfilled in general. But it is essential that $G(.;h,v_0)$ maps to $Z_T$ so that we can apply $L^{-1}$. For this reason we also subtracted the mean value in the second component of $G(w;h,v_0)$. In fact, \eqref{eq_navstk_transf1}-\eqref{eq_navstk_transf6} and the formulation \eqref{eq_fpgl_stokes} are equivalent because of
\begin{Remark}\label{th_bem_kompbed_äquiv}\upshape
Let $2<q<3, u\in W^1_{q,0}(\Omega)^n$ and $h\in U_0$ with $U_0$ as in Remark \ref{th_hanzawa_mod_bem}. Then
\[
\diverg\,u=g(h)u-\frac{1}{|\Omega|}\int_\Omega g(h)u\,dx\quad\Leftrightarrow\quad\diverg\,u=g(h)u.
\]
\end{Remark}
\textit{Proof.} The second equation is equivalent to $\diverg_h u=0$. For $u\in W^1_{q,0}(\Omega)^n$ holds $\int_\Omega\diverg\,u\,dx=0$ and this shows the reverse direction. The other part follows as in A. and Wilke \cite{Abels_Mullins_Sekerka}, p. 51.\hfill$\square$

The crucial step for showing Theorem \ref{th_navstk_exsatz} is to prove the following lemma which states the properties of the right hand side $G$.  
\begin{Lemma}\label{th_G}
Let $p>n+2$ and $2<q<3$ with $1+\frac{n+2}{p}>\frac{n+2}{q}$ as well as $0<T\leq 1$. Then for all $w\in Y_{T},\,h\in V_T$ and $v_0\in Y_{\gamma,h(0)}\cap L^q_\sigma(\Omega)$ we have that $G(w;h,v_0)$ is contained in $\tilde{Z}_T$ and for $w\in Y_{T,h(0),v_0}$ it holds $G(w;h,v_0)\in Z_T$. If $\|v_0\|_{Y_{\gamma,h(0)}}$ is bounded by a constant $R>0$, then $\|G(0;h,v_0)\|_{\tilde{Z}_T}\leq C_R$.

Moreover, let $R_1,R_2>0$. Then for $0<\varepsilon<R_0$ and $h_j\in V_{T,\varepsilon}$ with $h_0^j:=h_j(0)$, $v_0^j\in Y_{\gamma,h_j(0)}\cap L^q_\sigma(\Omega)$ and $\tilde{v}_0^j:=v_0^j\circ\tildeTheta_{h_j(0)}$ as well as $w_j=(u_j,\tilde{q}_j)\in Y_T$ with $u_0^j:=u_j(0)$, $\|u_0^j\|_{Y_\gamma}\leq R_1$ and $\|w_j\|_{Y_T}\leq R_2$ for $j=1,2$ the following estimate holds:
\begin{align*}\begin{split}
\|G(w_1;h_1,v_0^1)-G(w_2;h_2,v_0^2)\|_{\tilde{Z}_T}\leq C_{R_1,R_2}(T^\delta+\varepsilon)(\|w_1-w_2\|_{Y_T}+\|u_0^1-u_0^2\|_{Y_\gamma})+\\
+C_{R_1,R_2}(\|h_1-h_2\|_{\E_1(T)}+\|h_0^1-h_0^2\|_{X_\gamma})+\|\tilde{v}_0^1-\tilde{v}_0^2\|_{Y_\gamma},
\end{split}
\end{align*}
where $\delta:=\min\{\frac{\tau}{2},\frac{1}{q}-\frac{1}{p}\}>0$ and $\tau>0$ is as in Lemma \ref{th_A}.
\end{Lemma}

Also as preparation for the proof of Lemma \ref{th_G}, in the following we discuss where the conditions for the exponent $q$ come from:
\begin{Remark}\upshape\label{th_bem_regul_spur_u}
First of all, we have the restriction $2\leq q<3$ because of Theorem \ref{th_stokes_lin} below. Theorem \ref{th_navstk_exsatz} gives us the velocity field $u(h)\in Y_T^1$ in dependence of the height function $h$. Later we will insert it in equation \eqref{eq_transf_glg7} which should hold in $L^p(0,T;X_0)$. In order to show local well-posedness afterwards, we need enough regularity for $\tr_\Sigma\,u$. This means
\begin{align}\label{eq_regul_spur_u}
\tr_\Sigma\,u\in L^s(0,T;X_0)^n\quad\text{ for some }s>p.
\end{align}
Since $Y_\gamma=\{v\in L^p(\Omega):v|_{\Omega^\pm}\in W^{2-\frac{2}{q}}_q(\Omega^\pm)\}^n\cap W^1_{q,0}(\Omega)^n$, we get from Lemma \ref{th_einb_E_1T_mit_W2p_W1p0} below that
\begin{align}\label{eq_einb_Y_T1}
Y_T^1\hookrightarrow C^0([0,T];Y_\gamma)\cap L^q(0,T;W^2_q(\Omega_0)^n)
\end{align}
and the embedding constant is bounded independent of $T>0$ if we add the term $\|u(0)\|_{Y_\gamma}$ in the $Y_T^1$-norm. Now, for $\theta\in(0,1)$ we use the reiteration theorem in Lunardi \cite{Lunardi_Interpolation}, Corollary 1.24 to the result
\[
(W^{2-\frac{2}{q}}_q(\Omega^\pm),W^2_q(\Omega^\pm))_{\theta,q}=(L^q(\Omega^\pm),W^2_q(\Omega^\pm))_{(1-\theta)(1-\frac{1}{q})+\theta,q}=B_{q,q}^{2-2(1-\theta)\frac{1}{q}}(\Omega^\pm).
\]
Therefore Lemma \ref{th_intungl_einb_L_B} yields $Y_T^1|_{\Omega^\pm}\hookrightarrow L^s(0,T;B_{q,q}^{2-2(1-\theta)\frac{1}{q}}(\Omega^\pm)^n)$ with $s=\frac{q}{\theta}$. To achieve \eqref{eq_regul_spur_u} we want to choose $\theta\in(0,1)$ such that $\frac{q}{\theta}=s>p$ and that the Besov-space embeds into $W^1_p(\Omega^\pm)^n$. Because of $q\leq p$ the latter is fulfilled for 
\[
2-(1-\theta)\frac{2}{q}-\frac{n}{q}> 1-\frac{n}{p}\quad\Leftrightarrow\quad\theta> -\frac{q}{2}-\frac{nq}{2p}+\frac{n+2}{2}.
\]
Hence \eqref{eq_regul_spur_u} holds by the trace theorem if
\[
\frac{q}{p}>-\frac{q}{2}-\frac{nq}{2p}+\frac{n+2}{2}\quad\Leftrightarrow\quad \frac{q}{2}>-q\frac{n+2}{2p}+\frac{n+2}{2}\quad\Leftrightarrow\quad 1>\frac{n+2}{q}-\frac{n+2}{p}.
\]
Since $p>n+2$ necessary conditions on $q$ are $q>2$, if $n=2$, and $q>\frac{5}{2}$, if $n=3$. So we directly wrote $q>2$ instead of $q\geq 2$. On the other hand the conditions are fulfilled for many $p,q$.

Another point is that we have to estimate $u\cdot\nabla_hu$ in the proof of Lemma \ref{th_G} and therefore want to use the embedding $Y_T^1\hookrightarrow L^r(0,T;C^0(\overline{\Omega}))^n$ for some $r>q$. Embeddings over $\Omega^\pm$ yield $Y_\gamma\hookrightarrow C^0(\overline{\Omega})^n$ for $q>2$, if $n=2$ and $q>\frac{5}{2}$, if $n=3$, since then $2-\frac{2+n}{q}>0$ holds. Hence from \eqref{eq_einb_Y_T1} we even obtain $Y_T^1\hookrightarrow C^0([0,T];C^0(\overline{\Omega})^n)$ and the embedding constant is bounded independent of $T>0$ if we add $\|u(0)\|_{Y_\gamma}$ in the $Y_T^1$-norm.

Altogether, for simplicity we restrict ourselves in Theorem \ref{th_navstk_exsatz} and Lemma \ref{th_G} to exponents $p>n+2$ and $2<q<3$ with $1+\frac{n+2}{p}>\frac{n+2}{q}$. 
\end{Remark}
\textit{Proof of Lemma \ref{th_G}.} For simplicity we consider $0<T\leq1$. Properties of $\tildeTheta_h$ and $A(h)$ were shown in Lemma \ref{th_A}. The mean value theorem implies for all $h\in V_{T,\varepsilon}$ and $t\in[0,T]$
\begin{align}\begin{split}\label{eq_A_h_Absch}
\|A(h)|_t-I\|_{C^1(\overline{\Omega})^{n\times n}}&\leq \|A(h)|_t-A(h)|_0\|_{C^1(\overline{\Omega})^{n\times n}} + \|A(h)|_0-I\|_{C^1(\overline{\Omega})^{n\times n}}\leq\\
&\leq CT^{\tau}+c\|h(0)\|_{X_\gamma}\leq C (T^\tau+\varepsilon).\end{split}
\end{align}
Now let $w,h,v_0$ and $(u_j,\tilde{q}_j)=w_j, u_0^j$ as well as $h_j,h_0^j$ and $v_0^j,\tilde{v}_0^j$ for $j=1,2$ as in the lemma. We begin with the first component of $G$ which we have to estimate in $\tilde{Z}_T^1=L^q(0,T;L^q(\Omega)^n)$. Lemma \ref{th_A} implies $(\nabla-\nabla_h)\tilde{q}\in \tilde{Z}_T^1$ and \eqref{eq_A_h_Absch} yields
\begin{align*}\begin{split}
&\|(\nabla-\nabla_{h_1})\tilde{q}_1-((\nabla-\nabla_{h_2})\tilde{q}_2)\|_{\tilde{Z}_T^1}\leq \|(\nabla-\nabla_{h_1})(\tilde{q}_1-\tilde{q}_2)+(\nabla_{h_2}-\nabla_{h_1})\tilde{q}_2\|_{\tilde{Z}_T^1}\leq\\
&\leq C(T^\tau+\varepsilon)\|\tilde{q}_1-\tilde{q}_2\|_{Y_T^2}+C(\|h_1-h_2\|_{\E_1(T)}+\|h_0^1-h_0^2\|_{X_\gamma})R_2.
\end{split}\end{align*}
From Lemma \ref{th_A} we obtain $\mu^\pm(\diverg_h\nabla_hu-\Delta u)\in \tilde{Z}_T^1$ and null additions imply
\begin{align*}
\|&\mu^\pm(\diverg_{h_1}\nabla_{h_1}u_1-\Delta u_1)-\mu^\pm(\diverg_{h_2}\nabla_{h_2}u_2-\Delta u_2)\|_{\tilde{Z}_T^1}\leq\\
\leq c&(\|h_1-h_2\|_{\E_1(T)}+\|h_0^1-h_0^2\|_{X_\gamma})\|u_1\|_{Y_T^1}+c\|A(h_2)-I\|_{C^0([0,T];C^1(\overline{\Omega})^{n\times n})}\|u_1-u_2\|_{Y_T^1}.
\end{align*}
Together with this estimate we infer from \eqref{eq_A_h_Absch}
\begin{align*}
\|& a^\pm(h_1,D_x)w_1-a^\pm(h_2,D_x)w_2\|_{\tilde{Z}_T^1}\leq\\ 
&\leq C(T^\tau+\varepsilon)\|w_1-w_2\|_{Y_T}+ C_{R_2}(\|h_1-h_2\|_{\E_1(T)}+\|h_0^1-h_0^2\|_{X_\gamma}).
\end{align*}
Furthermore, Remark \ref{th_bem_regul_spur_u} yields $Y_T^1\hookrightarrow C^0([0,T];Y_\gamma)\hookrightarrow C^0([0,T];C^0(\overline{\Omega}))^n$ and the embedding constant is bounded independent of $T>0$ if we add $\|u(0)\|_{Y_\gamma}$ in the $Y_T^1$-norm. Hence we get $u\cdot\nabla_hv\in \tilde{Z}_T^1$ for $u,v\in Y_T^1$ because of Lemma \ref{th_A} and null additions imply
\begin{align*}
\|& u_1\cdot\nabla_{h_1}u_1-u_2\cdot\nabla_{h_2}u_2\|_{\tilde{Z}_T^1}\leq\\
&\leq C_{R_1,R_2}T^{\frac{1}{q}}(\|u_1-u_2\|_{Y_T^1}+\|u_0^1-u_0^2\|_{Y_\gamma})+C_{R_1,R_2}T^\frac{1}{q}(\|h_1-h_2\|_{\E_1(T)}+\|h_0^1-h_0^2\|_{X_\gamma}).
\end{align*}
Additionally we obtain $\nabla_hu\,\partial_t\tildeTheta_h\in\tilde{Z}_T^1$ from Lemma \ref{th_A}, the embedding for $Y_T^1$ mentioned above and $W^1_p(\Omega)\hookrightarrow C^0(\overline{\Omega})$ as $p>n+2$. Moreover, null additions yield
\begin{align*}
\|& \nabla_{h_1}u_1\,\partial_t\tildeTheta_{h_1}-\nabla_{h_2}u_2\,\partial_t\tildeTheta_{h_2}\|_{\tilde{Z}_T^1}\leq\\
&\leq C_{R_1,R_2}T^{\frac{1}{q}-\frac{1}{p}}(\|h_1-h_2\|_{\E_1(T)}+\|h_0^1-h_0^2\|_{X_\gamma})+ cT^{\frac{1}{q}-\frac{1}{p}}(\|u_1-u_2\|_{Y_T^1}+\|u_0^1-u_0^2\|_{Y_\gamma}).
\end{align*}

Now we consider the second component of $G$. First of all, Lemma \ref{th_A} implies that $g(h)u=\Tr(\nabla u-\nabla_hu)=(I-A(h)):\nabla u\in L^q(0,T;W^1_q(\Omega_0))$ and \eqref{eq_A_h_Absch} yields with a null addition
\[
\|g(h_1)u_1-g(h_2)u_2\|_{L^q(0,T;W^1_q(\Omega_0))}\leq c(T^\tau +\varepsilon)\|u_1-u_2\|_{Y_T^1}+c_{R_2}(\|h_1-h_2\|_{\E_1(T)}+\|h_0^1-h_0^2\|_{X_\gamma}).
\]
Therefore this also follows for $G(w;h,v_0)_2$ in $L^q(0,T;W^1_{q,(0)}(\Omega_0))$. Apart from that we have
$G(w;h,v_0)_2=g(h)u$ in $L^q(0,T;W^{-1}_{q,(0)}(\Omega))$ since constants drop out when integrating with $\phi\in W^1_{q',(0)}(\Omega)$. In this space we need to show the existence of a weak derivative. To this end let $\eta\in C_0^\infty(0,T)$. We approximate $u$ and $A(h)$ suitably. 
Using convolution one can show that there are $(u^m)_{m\in\N}$ in $C^\infty([0,T];W^2_q(\Omega_0)\cap W^1_{q,0}(\Omega))^n$ and $(A^m)_{m\in\N}$ in $C^\infty([0,T];C^1(\overline{\Omega})^{n\times n})$
with $u^m\rightarrow u$ in $Y_T^1$ and $A^m\rightarrow A(h)$ in the space $C^0([0,T];C^1(\overline{\Omega})^{n\times n})\cap W^1_p(0,T;L^p(\Omega)^{n\times n})$ for $m\rightarrow\infty$. Now the product rule yields
\[
C^\infty([0,T];W^1_q(\Omega_0))\ni g^m:=(I-A^m):\nabla u^m\overset{m\rightarrow\infty}{\longrightarrow} g(h)u\quad\text{ in }L^q(0,T;W^1_q(\Omega_0)).
\] 
Choose $\phi\in W^1_{q',(0)}(\Omega)$ arbitrary. Then for $m\in\N$ we obtain
\begin{align}\label{eq_bew_G_2}
-\int_0^T\partial_t\eta(t)\int_\Omega g^m(t)\phi\,dx\,dt=\int_0^T\eta(t)\int_\Omega \partial_t g^m(t)\phi\,dx\,dt.
\end{align}
Here $\partial_t g^m=-\partial_t A^m:\nabla u^m+(I-A^m):\nabla\partial_t u^m$ and therefore
\begin{align}\label{eq_bew_G_3}
\int_\Omega\partial_t g^m\phi\,dx=\int_\Omega-\partial_t A^m:\nabla u^m\,\phi+\partial_tu^m\cdot\diverg A^m\,\phi-((I-A^m)\partial_tu^m)\cdot\nabla\phi\,dx,
\end{align}
where we used integration by parts for the second term. In this equation we want to pass to the limit. Lemma \ref{th_einb_E_1T_mit_W2p_W1p0} implies $Y_T^1\hookrightarrow C^0([0,T];Y_\gamma)$ and the embedding constant is bounded independent of $T>0$ if we add $\|u(0)\|_{Y_\gamma}$ in the $Y_T^1$-norm. This yields
\[
\partial_tA^m:\nabla u^m\overset{m\rightarrow\infty}{\longrightarrow}\partial_tA(h):\nabla u\quad\text{ in }L^p(0,T;L^{r_1}(\Omega))\quad\text{ with }\frac{1}{r_1}:=\frac{1}{p}+\frac{1}{q}\in(0,1).
\] 
Furthermore we have $\partial_tu^m\cdot\diverg A^m\overset{m\rightarrow\infty}{\longrightarrow}\partial_tu\cdot \diverg\,A$ and $(I-A^m)\partial_tu^m\overset{m\rightarrow\infty}{\longrightarrow}(I-A(h))\partial_t u$ in $L^q(0,T;L^q(\Omega))$ (vector-valued for the second one). Now it holds $W^1_{q'}(\Omega)\hookrightarrow L^{r_2}(\Omega)$ for $n=2,3$ with $\frac{1}{r_2}:=\frac{1}{q'}-\frac{1}{n}\in(0,1)$ and because of
\[
\frac{1}{r_1}+\frac{1}{r_2}=1+\frac{1}{p}-\frac{1}{n}\in(0,1)\quad\text{ for }p>n+2, n=2,3
\]
the mapping $L^{r_1}(\Omega)\rightarrow W^{-1}_{q,(0)}(\Omega): g\mapsto \langle g,.\rangle$ is bounded and linear. This is also true for
\[
L^q(\Omega)\rightarrow W^{-1}_{q,(0)}(\Omega):g\mapsto\langle g,.\rangle\quad\text{ and }\quad L^q(\Omega)^n\rightarrow W^{-1}_{q,(0)}(\Omega): g\mapsto [\phi\mapsto\int_\Omega g\cdot\nabla\phi\,dx].
\] 
Since $\langle.,\phi\rangle:W^{-1}_{q,(0)}(\Omega)\rightarrow\R$ is bounded and linear for $\phi\in W^1_{q',(0)}(\Omega)$ we can pass to the limit in \eqref{eq_bew_G_2}-\eqref{eq_bew_G_3} for fixed $\phi$ and then we put outside the duality product. Hence $g(h)u\in W^1_q(0,T;W^{-1}_{q,(0)}(\Omega))$ and the derivative is given by
\begin{align}\label{eq_bew_G_4}
\langle\partial_t(g(h)u),\phi\rangle=\int_\Omega(-\partial_t A(h):\nabla u+\partial_tu\cdot\diverg A(h))\phi-((I-A(h))\partial_tu)\cdot\nabla\phi\,dx.
\end{align}
For the difference of the derivatives the above observations together with null additions and $\diverg\,I=0$ yield
\begin{align*}
\|&\partial_t(g(h_1)u_1)-\partial_t(g(h_2)u_2)\|_{L^q(0,T;W^{-1}_{q,(0)}(\Omega))}\leq\\
&\leq CT^{\frac{1}{q}-\frac{1}{p}}(C_{R_1,R_2}(\|h_1-h_2\|_{\E_1(T)}+\|h_0^1-h_0^2\|_{X_\gamma})+\|u_1-u_2\|_{Y_T^1}+\|u_0^1-u_0^2\|_{Y_\gamma})+\\
&+ C_{R_2}(\|h_1-h_2\|_{\E_1(T)}+\|h_0^1-h_0^2\|_{X_\gamma})+C(T^\tau+\varepsilon)\|u_1-u_2\|_{Y_T^1}.
\end{align*}
Now we show $\tr_\Sigma(G(w;h,v_0)_2|_{\Omega^+})\in\tilde{Z}_T^3$ and the related difference estimate. 
To this end let $\textup{ext}_T\in\Lc(W^{2,1}_q(\Omega^+_T),W^{2,1}_q(\Omega^+\times\R_+))$ for $T>0$ be extension operators in time such that
\[
\|\textup{ext}_T(u)\|_{W^{2,1}_q(\Omega^+\times\R_+)}\leq M(\|u\|_{W^{2,1}_q(\Omega^+_T)}+\|u(0)\|_{W^{2-\frac{2}{q}}_q(\Omega^+)})
\]
holds for some $M>0$. To construct such operators one can e.g. first extend suitably from $\Omega$ to $\R^n$ and then proceed as in Amann \cite{Amann_MaxReg_Quasilin}, Lemma 7.2 utilizing that the Laplace-operator has maximal $L^p$-regularity, cf. Prüss and Simonett \cite{PSimonett}, Theorem 6.1.8.

We can extend $\textup{ext}_Tu$ by reflection to a $\tilde{u}$ defined on $\Omega^+\times\R$. From Grubb \cite{Grubb}, (A.8)-(A.12) we know that $\nabla\tilde{u}\in H^{1,\frac{1}{2}}_q(\Omega^+\times\R)^{n\times n}\hookrightarrow W^{1,\frac{1}{2}}_q(\Omega^+\times\R)^{n\times n}$, where the latter embedding holds because of $q\in(2,3)$. Moreover, we have the estimate
\[
\|\nabla \tilde{u}\|_{W^{1,\frac{1}{2}}_{q}(\Omega^+\times\R)^{n\times n}}\leq C\|\tilde{u}\|_{W^{2,1}_{q}(\Omega^+\times\R)^n}\leq \tilde{C}(\|u\|_{Y_T^1}+\|u(0)\|_{Y_\gamma}).
\]
We decompose $A(h)=\tilde{A}(h)+A(h(0))$ with $\tilde{A}(h):=A(h)-A(h(0))$ and extend $\tilde{A}(h)$ by reflection at $T$ to $[0,2T]$ and then by $0$ to all of $\R_+$. Then $\tilde{A}(h)$ is contained in $C^{\frac{1}{2}+\tau}(\R_+;C^0(\overline{\Omega})^{n\times n})\cap C^\tau(\R_+;C^1(\overline{\Omega})^{n\times n})$ and the norm is estimated by the respective one of $A(h)-A(h(0))$. Since $\tilde{A}(h)=0$ on $(2T,\infty)$, we can apply Lemma \ref{th_hölder_prod_sob_bes} to $I=(0,2T+1)$ and $\tilde{A}(h)|_{\Omega^+}$. Afterwards we extend by $0$ and obtain
\[
(I-A(h(0))-\tilde{A}(h))|_{\Omega^+}:\nabla\tilde{u}\in W_q^{1,\frac{1}{2}}(\Omega^+\times\R_+).
\]
This then also holds for the mean value and e.g. Denk, Hieber and Prüss \cite{DenkHieberPruess_LpLq}, Lemma 3.5 implies
$\tr_\Sigma(G(w;h,v_0)_2|_{\Omega^+})\in \tilde{Z}_T^3$. Furthermore, a null addition yields
\begin{align*}
&\|\tr_\Sigma (G(w_1;h_1,v_0^1)_2|_{\Omega^+})-\tr_\Sigma(G(w_2;h_2,v_0^2)_2|_{\Omega^+})\|_{\tilde{Z}_T^3}\leq\\
&\leq C(\|I-A(h_0^1)\|_{C^1(\overline{\Omega})^{n\times n}}+\|\tilde{A}(h_1)\|_{C^{\frac{1+\tau}{2}}(\R_+;C^0(\overline{\Omega}))\cap C^0(\R_+;C^1(\overline{\Omega}))^{n\times n}})\cdot\\
& \cdot\|\tilde{u}_1-\tilde{u}_2\|_{W^{2,1}_{q}(\Omega^+\times\R_+)^n}+\|A(h_1)-A(h_2)\|_{C^{\frac{1}{2}+\tau}([0,T];C^0(\overline{\Omega}))\cap C^0([0,T];C^1(\overline{\Omega}))^{n\times n}}C_{R_1,R_2}\leq\\
&\leq C(T^{\tau/2}+\varepsilon)(\|u_1-u_2\|_{Y_T^1}+\|u_0^1-u_0^2\|_{Y_\gamma})+ C_{R_1,R_2}(\|h_1-h_2\|_{\E_1(T)}+\|h_0^1-h_0^2\|_{X_\gamma}),
\end{align*}
where we used Lemma \ref{th_A}, \eqref{eq_A_h_Absch} and $T\leq 1$.

Next we study the third component of $G$. Properties of $\nu_.$ were shown in Corollary \ref{th_A_kor}. As in \eqref{eq_A_h_Absch} we obtain $\|\nu_\Sigma-\nu_h\|_{C^0([0,T];C^1(\Sigma)^n)}\leq C(T^\tau+\varepsilon)$. Similar as before we get
$\lsprung 2\mu^\pm\textup{sym}(\nabla u)\rsprung ,\,\lsprung 2\mu^\pm\text{sym}(\nabla u-\nabla_h u)\rsprung\in(\tilde{Z}_T^3)^{n\times n}$
and the following estimates hold:
\begin{align}\begin{split}\label{eq_bew_G_11}
\|\lsprung 2\mu^\pm\textup{sym}(\nabla u)\rsprung\|_{(\tilde{Z}_T^3)^{n\times n}}& \leq C(\|u\|_{Y_T^1}+\|u(0)\|_{Y_\gamma}),\\
\|\lsprung 2\mu^\pm\text{sym}(\nabla u-\nabla_{h_1}u)\rsprung\|_{(\tilde{Z}_T^3)^{n\times n}}& \leq C(T^{\tau/2}+\varepsilon)(\|u\|_{Y_T^1}+\|u(0)\|_{Y_\gamma}),\\
\|\lsprung 2\mu^\pm\text{sym}(\nabla_{h_1}u-\nabla_{h_2}u)\rsprung\|_{(\tilde{Z}_T^3)^{n\times n}}& \leq C_{R_1,R_2}(\|h_1-h_2\|_{\E_1(T)}+\|h_0^1-h_0^2\|_{X_\gamma}).\end{split}
\end{align}
Since by definition $\lsprung\tilde{q}\rsprung\in \tilde{Z}_T^3$ and a product on $C^1(\Sigma)\times W_q^{1-\frac{1}{q}}(\Sigma)\rightarrow W_q^{1-\frac{1}{q}}(\Sigma)$ is given by pointwise multiplication, Lemma \ref{th_hölder_prod_sob_bes} implies 
\[
t^\pm(h;D_x)w=\lsprung 2\mu^\pm\textup{sym}(\nabla u)-\tilde{q}I\rsprung\,(\nu_\Sigma-\nu_h)+\lsprung 2\mu^\pm\text{sym}(\nabla u-\nabla_hu)\rsprung\,\nu_h\in(\tilde{Z}_T^3)^n.
\]
Corollary \ref{th_A_kor} yields $[\nu_h]_{C^{\frac{1+\tau}{2}}([0,T];C^0(\Sigma))^n}\leq CT^{\tau/2}$ and with null additions we infer
\begin{align*}
\|& t^\pm(h_1;D_x)w_1-t^\pm(h_2;D_x)w_2\|_{(\tilde{Z}_T^3)^n}\leq\\
&\leq C(T^{\tau/2}+\varepsilon)(\|w_1-w_2\|_{Y_T}+\|u_0^1-u_0^2\|_{Y_\gamma})+C_{R_1,R_2}(\|h_1-h_2\|_{\E_1(T)}+\|h_0^1-h_0^2\|_{X_\gamma}).
\end{align*}
Invoking Lemma \ref{th_mittl_krg2} and Corollary \ref{th_A_kor} together with Lemma \ref{th_hölder_prod_sob_bes} and the product rule, we obtain $\sigma K(.)\nu_. \in BC^1(V_T;(\tilde{Z}_T^3)^n)$ where $\sigma K(.)\nu_.$ and the derivative is bounded independent of $0<T\leq 1$ if we add $\|h(0)\|_{X_\gamma}$ in the $\E_1(T)$-norm. In particular we get
\[
\|\sigma K(h_1)\nu_{h_1}-\sigma K(h_2)\nu_{h_2}\|_{(\tilde{Z}_T^3)^n}\leq C(\|h_1-h_2\|_{\E_1(T)}+\|h_0^1-h_0^2\|_{X_\gamma}).
\]

For the fourth component we use that $\tildeTheta_{h(0)}:\Omega\rightarrow\Omega$ is a $C^2$-diffeomorphism and the identity in a neighbourhood of $\partial\Omega$ by Lemma \ref{th_hanzawa_mod} and Remark \ref{th_hanzawa_mod_bem}. Moreover, because of Lemma \ref{th_A} we have that $|\det\,A(h(0))|$ and the $C^2(\overline{\Omega})$-norm of $\tildeTheta_{h(0)}$ are bounded by a constant independent of $h\in V_T$. Therefore $W^k_q(\tildeTheta_{h(0)}(\Omega^\pm))^n\rightarrow W^k_q(\Omega^\pm)^n:v\mapsto v\circ\tildeTheta_{h(0)}$ for $k=0,1,2$ is bounded, linear and the operator norm is bounded independent of $h\in V_T$. The latter follows from Adams and Fournier \cite{Adams}, Theorem 3.41 and its proof. By interpolation this also holds for $Y_{\gamma,h(0)}\rightarrow Y_\gamma:v\mapsto v\circ\tildeTheta_{h(0)}$. We obtain $G(w;h,v_0)_4=v_0\circ\tildeTheta_{h(0)}\in Y_\gamma$ and, if $\|v_0\|_{Y_{\gamma,h(0)}}\leq R$ for some $R>0$, we have $\|v_0\circ\tildeTheta_{h(0)}\|_{Y_\gamma}\leq C_R$. The estimate for the difference of $\tilde{v}_0^j$ is trivial.

Altogether we have proven $G(w;h,v_0)\in\tilde{Z}_T$ and the difference estimate in the lemma. Additionally $G(0;h,v_0)=(0,\,0,\,-\sigma K(h)\nu_h,\,v_0\circ\tildeTheta_{h(0)})^\top$ and the above implies $\|G(0;h,v_0)\|_{\tilde{Z}_T}\leq C_R$. It remains to show for $w\in Y_{T,h(0),v_0}$ the compatibility condition for $G(w;h,v_0)$ in $Z_T$ which is equivalent to $\diverg\,(u(0))=(g(h)u)|_{t=0}$ in $W^{-1}_{q,(0)}(\Omega)$ since constants drop out when integrating with $\phi\in W^1_{q',(0)}(\Omega)$. Lemma \ref{th_einb_E_1T_mit_W2p_W1p0} implies that $Y_T^1\hookrightarrow C^0([0,T];Y_\gamma)$ and therefore
\[
g(h)u=(I-A(h)):\nabla u\in C^0([0,T];L^q(\Omega)).
\] 
Since $L^q(\Omega)\rightarrow W^{-1}_{q,(0)}(\Omega):g\mapsto\langle g,.\rangle$ is bounded and linear, it suffices to show $\diverg\,u(0)=g(h(0))u(0)$ or, equivalently, $\diverg_{h(0)}u(0)=0$ in $L^q(\Omega)$. But by chain rule we know $(\diverg\,v_0)\circ\tildeTheta_{h(0)}=\diverg_{h(0)}(u(0))$. Since $v_0$ is divergence-free, the claim follows.\hfill$\square$

After this preparation we can show the existence result for the Navier-Stokes part.

\noindent
{\textit{Proof of Theorem \ref{th_navstk_exsatz}.}} Let $p,q$ be as in the theorem and $R>0$ be arbitrary. Then for $h\in V_T$ and $v_0\in Y_{\gamma,h(0)}\cap L^q_{\sigma}(\Omega)$ with $\|v_0\|_{Y_{\gamma,h(0)}}\leq R$ Lemma \ref{th_G} yields
\begin{align}\label{eq_bew_navst1}
\|v_0\circ\tildeTheta_{h(0)}\|_{Y_\gamma}\leq\|G(0;h,v_0)\|_{\tilde{Z}_T}\leq C_R.
\end{align} 
We set $C_L:=\sup_{0<T\leq 1}\|L^{-1}\|_{\Lc(Z_T,Y_T)}$ and choose $R_1:=C_R$ and $R_2:=2C_RC_L$ in Lemma \ref{th_G}. Here $C_L$ is well-defined because of Theorem \ref{th_stokes_lin}. Then \eqref{eq_bew_navst1} implies for $0<\varepsilon< R_0$ and $0<T\leq 1$ as well as $h\in V_{T,\varepsilon}$ and $w\in Y_{T,h(0),v_0}$ with $\|w\|_{Y_T}\leq R_2$ the estimate
\begin{align*}
\|L^{-1}G(w;h,v_0)\|_{Y_T}&\leq \|L^{-1}\|_{\Lc(Z_T,Y_T)}(\|G(0;h,v_0)\|_{\tilde{Z}_T}+\|G(w;h,v_0)-G(0;h,v_0)\|_{\tilde{Z}_T})\leq\\
&\leq C_L(C_R+C_{R_1,R_2}(T^\delta+\varepsilon)(\|w\|_{Y_T}+\|v_0\circ\tildeTheta_{h(0)}\|_{Y_\gamma}))
\end{align*}
where $\delta>0$. For $T_0=T_0(R),\varepsilon=\varepsilon(R)>0$ small we obtain from \eqref{eq_bew_navst1} for all $0<T\leq T_0$ 
\begin{align}\label{eq_bew_navst2}
\|L^{-1}G(w;h,v_0)\|_{Y_T}\leq \frac{1}{2}R_2+\tilde{C}_R(T^\delta+\varepsilon)\leq\frac{2}{3}R_2.
\end{align}
Lemma \ref{th_G} yields $\|L^{-1}G(w_1;h,v_0)-L^{-1}G(w_2;h,v_0)\|_{Y_T}\leq C_L C_{R_1,R_2}(T_0^\delta+\varepsilon)\|w_1-w_2\|_{Y_T}$ for all $w_j\in Y_{T,h(0),v_0}$ with $\|w_j\|_{Y_T}\leq R_2$ for $j=1,2$. We choose $\varepsilon=\varepsilon(R)>0$ and $T_0=T_0(R)>0$ small such that
\begin{align}\label{eq_bew_navst3}
C_L C_{R_1,R_2}(T_0^\delta+\varepsilon)\leq \frac{1}{2}
\end{align}
holds. Then $L^{-1}G(.;h,v_0):\overline{B_{Y_{T,h(0),v_0}}(0,R_2)}\rightarrow \overline{B_{Y_{T,h(0),v_0}}(0,R_2)}$ is well-defined for all $0<T\leq T_0$ and a strict contraction. Here $Y_{T,h(0),v_0}$ is a complete metric space as closed subset of a Banach space. Banach's fixed-point theorem yields a unique fixed-point $F_T(h,v_0)$ in $\overline{B_{Y_{T,h(0),v_0}}(0,R_2)}$ and it holds
\[
F_{T}(h,v_0)|_{[0,\tilde{T}]}=F_{\tilde{T}}(h|_{[0,\tilde{T}]},v_0)\quad\text{ for all }0<\tilde{T}\leq T
\]
since $h|_{[0,\tilde{T}]}\in V_{\tilde{T},\varepsilon}$ as well as $F_T(h,v_0)|_{[0,\tilde{T}]}\in \overline{B_{Y_{\tilde{T},h(0),v_0}}(0,R_2)}$ and $L^{-1}G(.;h,v_0)$ is compatible with restrictions in time on $[0,\tilde{T}]$ for $0<\tilde{T}\leq T$. The uniqueness in $Y_{T,h(0),v_0}$ can be shown with the compatibility for restrictions in time and \eqref{eq_bew_navst2}.

It remains to prove the Lipschitz dependence. To this end let $h_j,h_0^j,v_0^j$ and $\tilde{v}_0^j$ for $j=1,2$ be as in the theorem. Then Lemma \ref{th_G} implies
\begin{align*}
\|& F_T(h_1,v_0^1)-F_T(h_2,v_0^2)\|_{Y_T}=\|L^{-1}(G(F_T(h_1,v_0^1);h_1,v_0^1)-G(F_T(h_2,v_0^2);h_2,v_0^2))\|_{Y_T}\leq\\
&\leq \|L^{-1}\|_{\Lc(Z_T,Y_T)}\|G(F_T(h_1,v_0^1);h_1,v_0^1)-G(F_T(h_2,v_0^2);h_2,v_0^2)\|_{Z_T}\leq\\
&\leq C_LC_{R_1,R_2}(T^\delta+\varepsilon)(\|F_T(h_1,v_0^1)-F_T(h_2,v_0^2)\|_{Y_T}+\|\tilde{v}_0^1-\tilde{v}_0^2\|_{Y_\gamma})+\\
&+C_LC_{R_1,R_2}(\|h_1-h_2\|_{\E_1(T)}+\|h_0^1-h_0^2\|_{X_\gamma})+\|\tilde{v}_0^1-\tilde{v}_0^2\|_{Y_\gamma}.
\end{align*}
Now from \eqref{eq_bew_navst3} the claim follows.\hfill$\square$

As preparation for the proof of the local well-posedness for \eqref{eq_transf_glg1}-\eqref{eq_transf_glg8} we show that the term coming from the Navier-Stokes part in \eqref{eq_transf_glg7} is sufficiently regular:
\begin{Corollary}\label{th_navst_kor}
Let $p,q$ and $\varepsilon, T_0$ for $R>0$ as in Theorem \ref{th_navstk_exsatz}. For $0<T\leq T_0, h\in V_{T,\varepsilon}$ and $v_0\in Y_{\gamma,h(0)}\cap L^q_\sigma(\Omega)$ with $\|v_0\|_{Y_{\gamma,h(0)}}\leq R$ for $j=1,2$ let
\[
G_T(h,v_0):=a_1(h)\cdot\tr_\Sigma(F_T(h,v_0)_1),
\]
where $a_1(h)$ is defined in Subsection \ref{sec_transf_glg} and $F_T(h,v_0)$ stems from Theorem \ref{th_navstk_exsatz}. Then for a $s>p$ it holds $G_T(h,v_0)\in L^s(0,T;X_0)$ and $\|G_T(h,v_0)\|_{L^s(0,T;X_0)}\leq C(R)$ as well as
\begin{align}\label{eq_G_kor_einschr}
G_T(h,v_0)|_{[0,\tilde{T}]}=G_{\tilde{T}}(h|_{[0,\tilde{T}]},v_0)\quad\text{ for all }0<\tilde{T}\leq T.
\end{align}
Moreover, for $h_j,h_0^j:=h_j(0),v_0^j$ and $\tilde{v}_0^j:=v_0^j\circ\tildeTheta_{h_0^j}, j=1,2$ as in Theorem \ref{th_navstk_exsatz} it holds
\[
\|G_T(h_1,v_0^1)-G_T(h_2,v_0^2)\|_{L^s(0,T;X_0)}\leq C_R(\|h_1-h_2\|_{\E_1(T)}+\|h_0^1-h_0^2\|_{X_\gamma}+\|\tilde{v}_0^1-\tilde{v}_0^2\|_{Y_\gamma}).
\]
\end{Corollary}
\textit{Proof.} From Remark \ref{th_bem_regul_spur_u} we know that $Y_T^1\rightarrow L^s(0,T;X_0)^n: u\mapsto \tr_\Sigma u$ is bounded, linear for some $s>p$ and the operator norm is bounded independent of $T>0$ if we add $\|u(0)\|_{Y_\gamma}$ in the $Y_T^1$-norm. Hence $\tr_\Sigma(F_T(h,v_0)_1)\in L^s(0,T;X_0)^n$ and from Theorem \ref{th_navstk_exsatz} and \eqref{eq_bew_navst1} we deduce $\|\tr_\Sigma(F_T(h,v_0)_1)\|_{L^s(0,T;X_0)^n}\leq C(R)$ as well as
\begin{align}\begin{split}\label{eq_bew_navstk_kor1}
\|&\tr_\Sigma(F_T(h_1,v_0^1)_1)-\tr_\Sigma(F_T(h_2,v_0^2)_1)\|_{L^s(0,T;X_0)^n}\leq\\
&\leq C_R(\|h_1-h_2\|_{\E_1(T)}+\|h_0^1-h_0^2\|_{X_\gamma}+\|\tilde{v}_0^1-\tilde{v}_0^2\|_{Y_\gamma}).\end{split}
\end{align}
Since pointwise multiplication is a product on $C^1(\Sigma)\times X_0\rightarrow X_0$, we obtain $G_T(h,v_0)\in L^s(0,T;X_0)$ with $\|G_T(h,v_0)\|_{L^s(0,T;X_0)}\leq \tilde{C}(R)$ by Corollary \ref{th_A_kor} and \eqref{eq_bew_navstk_kor1} together with a null addition yields the estimate in the lemma. Equation \eqref{eq_G_kor_einschr} follows from the respective one for $F_T$ in Theorem \ref{th_navstk_exsatz}.\hfill$\square$\\

\section{Local Well-Posedness}\label{sec_local_wellposed}
In this section we show the local well-posedness for the transformed Navier-Stokes/mean curvature flow system \eqref{eq_transf_glg1}-\eqref{eq_transf_glg8}. Let $p>n+2$ and $2<q<3$ with $1+\frac{n+2}{p}>\frac{n+2}{q}$. For $\varepsilon,T_0>0$ to $R>0$ as in Theorem \ref{th_navstk_exsatz} and $0<T\leq T_0$ as well as $h\in V_T, h_0\in X_\gamma$ with $\|h_0\|_{X_\gamma}\leq\varepsilon$ and $v_0\in Y_{\gamma,h(0)}\cap L^q_\sigma(\Omega)$ with $\|v_0\|_{Y_{\gamma,h(0)}}\leq R$ the system is equivalent to the abstract evolution equation
\begin{alignat}{2}
\partial_t h-G_T(h,v_0)&=a_2(h)\,K(h) &\qquad&\text{ on }\Sigma\times(0,T),\label{eq_evolgl_h1}\\
h|_{t=0}&=h_0 &\qquad&\text{ on }\Sigma,\label{eq_evolgl_h2}
\end{alignat}
where $G_T$ is defined in Corollary \ref{th_navst_kor}. Here \enquote{equivalent} means that $(u,\tilde{q},h)\in Y_T\times V_T$ solves \eqref{eq_transf_glg1}-\eqref{eq_transf_glg8} if and only if $h$ fulfils \eqref{eq_evolgl_h1}-\eqref{eq_evolgl_h2} and it holds $(u,\tilde{q})=F_T(h,v_0)$.

In order to show local well-posedness for \eqref{eq_evolgl_h1}-\eqref{eq_evolgl_h2}, we use that by Lemma \ref{th_mittl_krg} the transformed mean curvature $K(h)$ has a quasilinear structure with respect to the height function $h$ and for $h=0$ the principal part is given by the Laplace-Beltrami operator $\Delta_\Sigma$. Here $\Delta_\Sigma:X_1\rightarrow X_0$ has maximal $L^p$-regularity on finite intervals because of Theorem 6.4.3 in Prüss and Simonett \cite{PSimonett}. To give a simple proof for the local well-posedness, we proceed as in the case of quasilinear parabolic equations, compare e.g. Köhne, Prüss and Wilke \cite{Koehne_Pruess_Wilke}, Theorem 2.1 (for $\mu=1$ there). We especially cannot apply this result directly since $G_T$ is a non-local operator. That means $G_T(h,v_0)$ at point $x_0\in\Sigma$ and time $t_0\in[0,T]$ in general depends on values $h(x,t)$ with $(x,t)\in\Sigma\times[0,T]$ outside a neighbourhood of $(x_0,t_0)$. Nevertheless we can essentially adapt the proof of \cite{Koehne_Pruess_Wilke} (for $\mu=1$) since by Corollary \ref{th_navst_kor} we have that $G_T$ is sufficiently regular and compatible with restrictions in time on $[0,\tilde{T}]$ for $0<\tilde{T}\leq T$. Intuitively the velocity and pressure at time $T$ depend on the evolution of the entire interface from the beginning up to time $T$.

The result for the local well-posedness of the transformed Navier-Stokes/mean curvature flow system \eqref{eq_transf_glg1}-\eqref{eq_transf_glg8} reads as follows:
\begin{Theorem}\label{th_locw}
Let $p>n+2$ and $2<q<3$ with $1+\frac{n+2}{p}>\frac{n+2}{q}$. For $R>0$ there exist $T_0=T_0(R),\varepsilon_0=\varepsilon_0(R)>0$ such that for all $0<T\leq T_0$ and $h_0\in X_\gamma$ with $\|h_0\|_{X_\gamma}\leq\varepsilon_0$ and $v_0\in Y_{\gamma,h_0}\cap L^q_\sigma(\Omega)$ with $\|v_0\|_{Y_{\gamma,h_0}}\leq R$ there is a unique solution $(u,\tilde{q},h)(h_0,v_0)\in Y_T\times V_T$ of \eqref{eq_transf_glg1}-\eqref{eq_transf_glg8} and
\[
\|(u,\tilde{q},h)(h_0^1,v_0^1)-(u,\tilde{q},h)(h_0^2,v_0^2)\|_{Y_T\times V_T}\leq C_R(\|h_0^1-h_0^2\|_{X_\gamma}+\|\tilde{v}_0^1-\tilde{v}_0^2\|_{Y_\gamma})
\]
holds for $h_0^j\in X_\gamma$ with $\|h_0^j\|_{X_\gamma}\leq\varepsilon_0$ and $v_0^j\in Y_{\gamma,h_j(0)}\cap L^q_\sigma(\Omega)$ with $ \|v_0^j\|_{Y_{\gamma,h_j(0)}}\leq R$ as well as $\tilde{v}_0^j:=v_0^j\circ\tildeTheta_{h_0^j}$ for $j=1,2$.
\end{Theorem}
\textit{Proof.} With $P,Q$ as in Lemma \ref{th_mittl_krg} and $U_0$ as in Remark \ref{th_hanzawa_mod_bem} we denote $\tilde{B}(h):=a_2(h)\,P(h)$ and $\tilde{F}(h):=a_2(h)\,Q(h)$ for all $h\in U_0$. Because of Theorem 6.4.3 in Prüss and Simonett \cite{PSimonett} there is a $w>0$ such that $\Delta_\Sigma-wI:X_1\rightarrow X_0$ has maximal $L^p$-regularity on $\R_+$. Therefore we consider $B:=\tilde{B}-wI$ and $ F:=\tilde{F}+wI$. Since $a_2(0)=1$ we have $B(0)=\Delta_\Sigma-wI$. Moreover, Corollary \ref{th_A_kor} yields $B\in C^1(U_0;\Lc(X_1,X_0))$ and $F\in C^1(U_0;X_0)$. The embedding $X_\gamma\hookrightarrow L^\infty(\Sigma)$ implies $\overline{B_{X_\gamma}(0,\delta_0)}\subseteq U_0$ for $0<\delta_0<R_0$ small. So if $\delta_0$ is sufficiently small, by the mean value theorem, there is a $L>0$ such that
\begin{align}\begin{split}\label{eq_bew_locw_lip}
\|B(h_1)v-B(h_2)v\|_{X_0}&\leq L\|h_1-h_2\|_{X_\gamma}\|v\|_{X_1},\\
\|F(h_1)-F(h_2)\|_{X_0}&\leq L\|h_1-h_2\|_{X_\gamma}
\end{split}\end{align}
holds for all $h_1,h_2\in \overline{B_{X_\gamma}(0,\delta_0)}$ and $v\in X_1$.

Now fix $R>0$ and let $T_0=T_0(R),\varepsilon=\varepsilon(R)>0$ be as in Theorem \ref{th_navstk_exsatz}. Then for $0<T\leq T_0$ and $h\in V_T$ as well as $h_0\in X_\gamma$ with $\|h_0\|_{X_\gamma}\leq\varepsilon$ and $v_0\in Y_{\gamma,h_0}\cap L^q_\sigma(\Omega)$ with $\|v_0\|_{Y_{\gamma,h_0}}\leq R$ the system \eqref{eq_transf_glg1}-\eqref{eq_transf_glg8} or \eqref{eq_evolgl_h1}-\eqref{eq_evolgl_h2}, respectively, are equivalent to
\begin{align}
\partial_th-B(h)h&=F(h)+G_T(h,v_0)\quad\text{ in }(0,T),\label{eq_evolgl_h1_abstr}\\
h(0)&=h_0\label{eq_evolgl_h2_abstr}.
\end{align}
Here $B(h)$ and $F(h)$ are defined in the sense of $V_T\subseteq C^0([0,T];U_0)$ by Remark \ref{th_hanzawa_mod_bem}. In particular, it holds $B(h)\in C^0([0,T];\Lc(X_1,X_0))$.

The idea for solving \eqref{eq_evolgl_h1_abstr}-\eqref{eq_evolgl_h2_abstr} is as follows: we fix $B(h)$ at $h=0$ and apply Banach's fixed-point theorem to $H(.;h_0,v_0)$ where $H(h;h_0,v_0):=f$ is the unique solution of
\begin{align}\label{eq_bew_locw_defH1}
\partial_tf-B(0)f&=B(h)h-B(0)h+F(h)+G_T(h,v_0)\quad\text{ in }(0,T),\\
f(0)&=h_0.\label{eq_bew_locw_defH2}
\end{align}
The unique solvability of \eqref{eq_bew_locw_defH1}-\eqref{eq_bew_locw_defH2} will follow from Prüss and Simonett \cite{PSimonett}, Theorem 3.5.5 if we extend the right hand side by zero to $\R_+$. To this end we need suitable sets for the height function. For $0<\delta<\min\{\varepsilon,\delta_0\}, 0<r<R_0-\delta_0$ and $0<T\leq T_0$ we set
\[
\B_{r,T,h_0}:=\{h\in\E_1(T)\,:\,h(0)=h_0, \|h\|_{\E_1(T)}\leq r\}\quad\text{ for }h_0\in \overline{B_{X_\gamma}(0,\delta)}.
\]
For $r,\delta>0$ small we show $\B_{r,T,h_0}\subseteq V_T\cap C^0([0,T];\overline{B_{X_\gamma}(0,\delta_0)})$. Therefore let $h\in\B_{r,T,h_0}$.

Because of $\|h\|_{\E_1(T)}+\|h_0\|_{X_\gamma}\leq r+\delta<R_0$ and our choice of $\delta_0$ it suffices to show $\|h(t)\|_{X_\gamma}\leq\delta_0$ for all $t\in[0,T]$ if $r,\delta>0$ are small independent of $h$. By Prüss and Simonett \cite{PSimonett}, Theorem 3.5.5 there is a unique solution $h_0^\ast\in\E_1(\R_+)$ of
\begin{align}\label{eq_bew_locw_hstern}
\partial_th_0^\ast-B(0) h_0^\ast=0,\quad h_0^\ast(0)=h_0.
\end{align}
Since $h(0)=h_0^\ast(0)$, we get from Theorem \ref{th_einb_aus_trace_method} with a null addition
\begin{align}\begin{split}\label{eq_bew_locw1}
\|h\|_{C^0([0,T];X_\gamma)}\leq C(\|h\|_{\E_1(T)}+\|h_0\|_{X_\gamma})\leq c(r+\delta).
\end{split}\end{align}
For $r,\delta>0$ small we obtain $\|h(t)\|_{X_\gamma}\leq\delta_0$ for all $t\in[0,T]$ and thus the above inclusion.

Now we can define $H(.;h_0,v_0)\equiv H_{r,T}(.;h_0,v_0):\B_{r,T,h_0}\rightarrow\E_1(T):h\mapsto f$ where $f$ is the unique solution to \eqref{eq_bew_locw_defH1}-\eqref{eq_bew_locw_defH2}. Then $h\in\B_{r,T,h_0}$ is a solution of \eqref{eq_evolgl_h1_abstr}-\eqref{eq_evolgl_h2_abstr} if and only if $h$ is a fixed-point of $H(.;h_0,v_0)$. To apply Banach's theorem we need the following:

For $r,\delta,T>0$ small, all $g_j\in \overline{B_{X_\gamma}(0,\delta)}$ and $h_j\in \B_{r,T,g_j}$ as well as $v_0^j\in Y_{\gamma,g_j}\cap L^q_\sigma(\Omega)$ with $\|v_0^j\|_{Y_{\gamma,g_j}}\leq R$ and $\tilde{v}_0^j:=v_0^j\circ\tildeTheta_{g_j}$ for $j=1,2$ we will show the following inequality:
\[
\|H(h_1;g_1,v_0^1)-H(h_2;g_2,v_0^2)\|_{\E_1(T)}\leq\frac{1}{2}\|h_1-h_2\|_{\E_1(T)}+c_R(\|g_1-g_2\|_{X_\gamma}+\|\tilde{v}_0^1-\tilde{v}_0^2\|_{Y_\gamma}).
\]
Moreover, for all $r>0$ there are $\delta(r),T(r)>0$ such that for all $0<\delta\leq\delta(r)$ and $0<T\leq T(r)$ the estimate $\|H(0;h_0,v_0)\|_{\E_1(T)}\leq \frac{r}{2}$ is valid, where $h_0\in \overline{B_{X_\gamma}(0,\delta)}$.

Now let us prove this. From Prüss and Simonett \cite{PSimonett}, Theorem 3.5.5 we know that there exists a $C>0$ independent of $T>0$ such that
\begin{align*}
&\|H(h_1;g_1,v_0^1)-H(h_2;g_2,v_0^2)\|_{\E_1(T)}\leq C(\|g_1-g_2\|_{X_\gamma}+\|F(h_1)-F(h_2)\|_{\E_0(T)}+\\
&+\|(B(h_1)-B(0))(h_1-h_2)+(B(h_1)-B(h_2))h_2+G_T(h_1,v_0^1)-G_T(h_2,v_0^2)\|_{\E_0(T)}).
\end{align*}
Here because of \eqref{eq_bew_locw_lip} it holds $\|F(h_1)-F(h_2)\|_{\E_0(T)}\leq LT^\frac{1}{p}\|h_1-h_2\|_{C^0([0,T];X_\gamma)}$ and with $g_j^\ast$ as in \eqref{eq_bew_locw_hstern} for $g_j, j=1,2$ instead of $h_0$ we conclude from Theorem \ref{th_einb_aus_trace_method}
\begin{align}\begin{split}\label{eq_bew_locw2}
\|h_1-h_2\|_{C^0([0,T];X_\gamma)}\leq c(\|h_1-h_2\|_{\E_1(T)}+\|g_1-g_2\|_{X_\gamma}).
\end{split}\end{align} 
Furthermore \eqref{eq_bew_locw_lip} together with \eqref{eq_bew_locw1} and \eqref{eq_bew_locw2} implies 
\begin{align*}
\|&(B(h_1)-B(0))(h_1-h_2)+(B(h_1)-B(h_2))h_2\|_{\E_0(T)}\leq\\ 
&\leq L(\|h_1\|_{C^0([0,T];X_\gamma)}\|h_1-h_2\|_{\E_1(T)}+\|h_1-h_2\|_{C^0([0,T];X_\gamma)}\|h_2\|_{\E_1(T)})\leq\\
&\leq LC(r+\delta)\|h_1-h_2\|_{\E_1(T)}+Lcr(\|h_1-h_2\|_{\E_1(T)}+\|g_1-g_2\|_{X_\gamma}).
\end{align*}
Additionally, for some $s>p$ Corollary \ref{th_navst_kor} and Hölder's inequality yield
\[
\|G_T(h_1,v_0^1)-G_T(h_2,v_0^2)\|_{\E_0(T)}\leq C_RT^{\frac{1}{p}-\frac{1}{s}}(\|h_1-h_2\|_{\E_1(T)}+\|g_1-g_2\|_{X_\gamma}+\|\tilde{v}_0^1-\tilde{v}_0^2\|_{Y_\gamma}).
\]
For $r,\delta,T>0$ small we obtain the claimed difference estimate. It remains to treat $\|H(0;h_0,v_0)\|_{\E_1(T)}$. Theorem 3.5.5 in Prüss and Simonett \cite{PSimonett} and Corollary \ref{th_navst_kor} imply
\begin{align*}
\|H(0;h_0,v_0)\|_{\E_1(T)}\leq C(\delta+\|F(0)\|_{\E_0(T)}+T^{\frac{1}{p}-\frac{1}{s}}C_R).
\end{align*}
So if $\delta(r),T(r)>0$ are appropriate, we get $\|H(0;h_0,v_0)\|_{\E_1(T)}\leq \frac{r}{2}$ for $0<\delta\leq\delta(r)$ and $0<T\leq T(r)$. Altogether we have proven everything claimed above.

After this preparation we can choose $r,\delta,T_0>0$ small such that everything before holds for $r,\delta,T_0$ and for $2r,\delta,T_0$. Then this also follows for $\tilde{r},\delta,T$ where $\tilde{r}\in\{r,2r\}$ and $0<T\leq T_0$. Now we choose $h_0=g_1=g_2\in \overline{B_{X_\gamma}(0,\delta)}$ and $v_0=v_0^1=v_0^2$ above and obtain that $H_{\tilde{r},T}(.;h_0,v_0)$ is a strict contraction on $\B_{\tilde{r},T,h_0}$ where the latter is a complete metric space as closed subset of a Banach space. Banach's theorem now yields a unique fixed-point $h^\ast_{\tilde{r},T}(h_0,v_0)\in\B_{\tilde{r},T,h_0}$ of $H_{\tilde{r},T}(.;h_0,v_0)$ in $\B_{\tilde{r},T,h_0}$ which is a solution of \eqref{eq_evolgl_h1_abstr}-\eqref{eq_evolgl_h2_abstr} and 
\[
(u^\ast,\tilde{q}^\ast,h^\ast)(h_0,v_0):=(F_T(h^\ast_{r,T}(h_0,v_0),v_0),h^\ast_{r,T}(h_0,v_0))\in Y_T\times V_T
\] 
solves the transformed Navier-Stokes/mean curvature flow system \eqref{eq_transf_glg1}-\eqref{eq_transf_glg8}. The Lipschitz-dependence of $(u^\ast,\tilde{q}^\ast,h^\ast)$ on $(h_0,\tilde{v}_0)$ where $\tilde{v}_0:=v_0\circ\tildeTheta_{h_0}$ is obtained as follows: For $h^\ast$ we use $H_{r,T}(h^\ast(h_0,v_0);h_0,v_0)=h^\ast(h_0,v_0)$ in the difference estimate proven before. Then for $(u^\ast,\tilde{q}^\ast)$ we utilize Theorem \ref{th_navstk_exsatz}.

It remains to show uniqueness in $Y_T\times V_T$ for fixed $h_0\in\overline{B_{X_\gamma}(0,\delta)}$ and $v_0$ as in the theorem. In the following we write $h^\ast_{\tilde{r},T}$ and $(u^\ast,\tilde{q}^\ast,h^\ast)$ instead of $h^\ast_{\tilde{r},T}(h_0,v_0)$ and $(u^\ast,\tilde{q}^\ast,h^\ast)(h_0,v_0)$, respectively. In particular it holds $h^\ast=h^\ast_{r,T}$. We need that the fixed-point $h^\ast_{\tilde{r},T}$ is the same for $\tilde{r}\in\{r,2r\}$ and remains the same up to restriction in time if we shrink $T$. Therefore let $0<\tilde{T}\leq T$ and $\tilde{r}\in\{r,2r\}$. For $h\in\B_{\tilde{r},T,h_0}$ we have $h|_{[0,\tilde{T}]}\in\B_{\tilde{r},\tilde{T},h_0}$ and because of Corollary \ref{th_navst_kor} it holds $G_T(h,v_0)|_{[0,\tilde{T}]}=G_{\tilde{T}}(h|_{[0,\tilde{T}]},v_0)$. Hence we obtain
\[
H_{\tilde{r},T}(h;h_0,v_0)|_{[0,\tilde{T}]}=H_{\tilde{r},\tilde{T}}(h|_{[0,\tilde{T}]};h_0,v_0).
\] 
Since $H_{r,\tilde{T}}(h;h_0,v_0)=H_{2r,\tilde{T}}(h;h_0,v_0)$ for all $h\in\B_{r,\tilde{T},h_0}\subseteq\B_{2r,\tilde{T},h_0}$ the uniqueness of the respective fixed-point yields
\[
h^\ast|_{[0,\tilde{T}]}=h^\ast_{r,\tilde{T}}=h^\ast_{2r,\tilde{T}}\quad\text{ for all }0<\tilde{T}\leq T.
\]
Now the uniqueness in $Y_T\times V_T$ can be shown by a typical contradiction argument.\hfill$\square$\\\appendix

\section{Banach-Space-Valued Functions}\label{sec_fctsp_vecvalued}

\subsection{Spaces of (Hölder-)Continuous Functions}\label{sec_hölderräume}
Let $X$ be a Banach space over $\K=\R$ or $\C$. For a closed interval $I\subseteq\R$ we denote by $C^0(I;X), C_b^0(I;X), BUC(I;X)$ and $C^{0,\alpha}(I;X)$ for $\alpha\in(0,1]$ the set of continuous; continuous and bounded; bounded and uniformly continuous; bounded and Hölder-continuous functions $f:I\rightarrow X$, respectively. In the latter case we write $[f]_{C^{0,\alpha}(I;X)}$ for the semi-norm and we set $C^{\alpha}(I;X):=C^{0,\alpha}(I;X)$ if $\alpha\in(0,1)$. Moreover, let $C^{0,0}(I;X):=C^0(I;X)$. If $I\subseteq\R$ is an open interval, then $C^\infty(\overline{I};X)$ is the set of all smooth $f:I\rightarrow X$ such that $f$ and all derivatives can be extended continuously to $\overline{I}$ and $C_0^\infty(I;X)$ denotes the set of all $f\in C^\infty(\overline{I};X)$ with compact support in $I$. For suitable product estimates we need:
\begin{Definition}\label{th_def_produkt}\upshape
Let $X,Y,Z$ be Banach spaces over $\K=\R$ or $\C$. Then $B:X\times Y\rightarrow Z$ is called \textit{product} if $B$ is bilinear and $\|B(f,g)\|_Z\leq C_0\|f\|_X\|g\|_Y$ holds for all $(f,g)\in X\times Y$ and some $C_0>0$.
\end{Definition}

With a null addition one can show that for any closed, bounded interval $I\subseteq\R$ and $\alpha\in[0,1]$ the product $B$ induces a product on the corresponding Hölder-spaces to the exponent $\alpha$ and in the estimate one can choose the same $C_0$. For estimates of nonlinear terms we use
\begin{Lemma}\label{th_hölderräume_nichtlin}
Let $I\subseteq\R$ a closed, bounded interval, $n,m\in\N$ and $\Omega\subseteq\R^n$ open and bounded, $U\subseteq\R^m$ open, $F:U\rightarrow\R$ be $C^1$ and $K\subseteq U$ compact and convex. For $\alpha\in[0,1]$ let
\[
M:=\{u\in C^{0,\alpha}(I;C^0(\overline{\Omega};\R^m)):u(t)(x)\in K\,\text{ for all }t\in I,x\in\overline{\Omega}\}
\]
and $\tilde{F}(u)(t)(x):=F(u(t)(x))$ for $u\in M$ and $t\in I,x\in\overline{\Omega}$. Then for $u\in M$ with  $\|u\|_{C^{0,\alpha}(I;C^0(\overline{\Omega})^m)}\leq R$ it holds $\tilde{F}(u)\in C^{0,\alpha}(I;C^0(\overline{\Omega}))$ with $\|\tilde{F}(u)\|_{C^{0,\alpha}(I;C^0(\overline{\Omega}))}\leq C(R)$ where $C(R)>0$ is independent of $I$ and $\alpha$.

Moreover, if $F$ is $C^2$, then for all $u,v\in M$ with norm estimated by $R$ there is a $C(R)>0$ independent of $I$ and $\alpha$ such that $\|\tilde{F}(u)-\tilde{F}(v)\|_{C^{0,\alpha}(I;C^0(\overline{\Omega}))}\leq C(R)\|u-v\|_{C^{0,\alpha}(I;C^0(\overline{\Omega})^m)}$.
\end{Lemma}
\textit{Proof.} 
For $u,v\in M$ and $t\in I,x\in\overline{\Omega}$ the mean value theorem implies
\[
\tilde{F}(u(t)(x))-\tilde{F}(v(t)(x))=\int_0^1 DF(su(t)(x)+(1-s)v(t)(x))\,ds\cdot (u(t)(x)-v(t)(x)).
\]
Now the first part directly follows from the compactness of $K$. If additionally $F$ is $C^2$ we apply the first part to $G:U\times U\rightarrow\R:(\tilde{u},\tilde{v})\mapsto \int_0^1 DF(s\tilde{u}+(1-s)\tilde{v})\,ds$ and $K\times K$ instead of $K$ using a product estimate.\hfill$\square$

In the following lemma let all functions have values in $\R$.
\begin{Lemma}\label{th_hölder_ableiten}
Let $I\subseteq\R$ be a closed, bounded interval, $\Omega\subseteq\R^n$ be a bounded, connected domain, $\alpha\in[0,1]$ and $k=0,1$. For $\beta\in\R$ and $c_0>0$ we consider $F(h):=h^\beta$ acting on $U_{k}:=\{ h\in C^k(\overline{\Omega}): h(x)>c_0\text{ for all }x\in\overline{\Omega}\}$ and on 
\[
V_{\alpha,k}:=\{h\in C^{0,\alpha}(I;C^k(\overline{\Omega})): h(t)\in U_k\text{ for all }t\in I\},
\] 
respectively. Then $F\in C^1(U_k;C^k(\overline{\Omega}))$ and $F\in C^1(V_{\alpha,k};C^{0,\alpha}(I;C^k(\overline{\Omega})))$ and for $R>0$ arbitrary and all $h\in U_k$ and $h\in V_{\alpha,k}$ with norm less or equal $R$, respectively, it holds that $F(h)$ and $F'(h)$ is bounded by a constant $C_{R,\beta,c_0}>0$ independent of $I$.
\end{Lemma}
\textit{Proof.} For $\beta\in\N_0$ the claim follows from the product rule. Therefore it is enough to consider the case $-\infty<\beta< 1$. By scaling we can also assume $c_0=1$. Now one uses that for $\beta<1$ and $x,y\geq 1$ it holds that $|x^\beta-y^\beta|\leq c_\beta|x-y|$ and $x^\beta\leq C_\beta(1+|x|)$. Moreover, for $k=1$ we have $\partial_{x_i}(F(h))=\beta\,h^{\beta-1}\partial_{x_i}h$ for $ i=1,...,n$. From this the continuity of $F$ on respective spaces and the estimate for $F(h)$ follows. In both cases our candidate for the derivative is $G(h)(\rho):=\beta h^{\beta-1}\rho$ on corresponding spaces. Replacing $\beta-1$ by $\beta$ above one can show that $G$ is well-defined, continuous and satisfies the desired estimate for $F'$ in the lemma. To verify the definition of the Fréchet-derivative one uses that for small $\rho$
\[
(F(h+\rho)-F(h)-G(h)\rho)=\int_0^1\int_0^1\beta(\beta-1)(h+uv\rho)^{\beta-2}u\,dv\,du\,\,\rho^2.
\]
If $\rho$ is bounded by a small $\varepsilon>0$ in the respective norm, Lemma \ref{th_hölderräume_nichtlin} yields for $k=0$ that the first part can be estimated by $c(h,\varepsilon)$ in the corresponding norm. Then product estimates yield the case $k=0$ and the case $k=1$ follows from this by replacing $\beta-1$ by $\beta$ and using the identity for $\partial_{x_i}(F(h))$.\hfill$\square$\\

\subsection{Lebesgue-, Sobolev- and Besov-Spaces}\label{sec_leb_sob_besov}
Let $I\subseteq\R$ be measurable, $X$ be a Banach space and $1\leq p\leq\infty$. Then $L^p(I;X)$ are the usual Bochner spaces and $W^k_p(I;X)$ for $I$ open and $k\in\N$ are the $X$-valued Sobolev spaces.

Moreover, we need vector-valued variants of some fractional order Besov- and Sobolev-spaces: Therefore let $s\in(0,1), 1\leq p\leq\infty$ and $I$ be an interval in $\R$. Then we define $B^s_{p,\infty}(I;X):=\{f\in L^p(I;X):\|f\|_{B^s_{p,\infty}(I;X)}<\infty\}$ where the norm is given by $\|f\|_{B^s_{p,\infty}(I;X)}:=\|f\|_{L^p(I;X)}+[f]_{B^s_{p,\infty}(I;X)}$ with
\[
[f]_{B^s_{p,\infty}(I;X)}:=\sup_{0<h\leq 1}\left\|\frac{f(.+h)-f}{h^s}\right\|_{L^p(I_h;X)}\quad\text{ and }I_h:=\{x\in I:x+h\in I\}.
\] 
For $s\in(0,1)$ and $1\leq p<\infty$ we set $W^s_p(I;X):=\{f\in L^p(I;X):\|f\|_{W^s_p(I;X)}<\infty\}$ where the norm is given by $\|f\|^p_{W^s_p(I;X)}:=\|f\|_{L^p(I;X)}^p+[f]^p_{W^s_p(I;X)}$ and
\[
[f]^p_{W^s_p(I;X)}:=\int_I\int_I\left(\frac{\|f(x)-f(y)\|}{|x-y|^s}\right)^p\frac{dx\,dy}{|x-y|}.
\]
All the above spaces are Banach spaces since $L^p(I;X)$ and $L^1(I\times I)$ are complete. The above definitions equal those in Simon \cite{Simon} up to equivalent norms. Moreover, for $s>0$ and $1<p<\infty$ the Bessel-potential spaces are denoted by $H^s_p(I;X)$.

Additionally, we need spaces with mixed regularity. Let $\Omega\subseteq\R^n$ be open, $r,s\in\N$, $1<p<\infty$ and $T>0$. Then we define $W_p^{r,s}(\Omega_T):=L^p(0,T;W_p^r(\Omega))\cap W^s_p(0,T;L^p(\Omega))$. If $T=\infty$ we write $\Omega\times\R_+$ instead of $\Omega_T$. Further anisotropic spaces are defined analogously. For properties of such spaces cf. the appendix in Grubb \cite{Grubb} and Denk, Hieber and Prüss \cite{DenkHieberPruess_LpLq}.

We need product estimates for Hölder spaces with suitable Sobolev- and Besov-spaces:
\begin{Lemma}\label{th_hölder_prod_sob_bes}
Let $X,Y,Z$ be Banach spaces over $\K=\R$ or $\C$ and $B:X\times Y\rightarrow Z$ be a product in sense of Definition \ref{th_def_produkt}. Moreover, let $I\subseteq\R$ be a finite interval, $s\in(0,1)$ and $1\leq p\leq\infty$. Then $C^s(\overline{I};X)\times B^s_{p,\infty}(I;Y)\rightarrow B^s_{p,\infty}(I;Z):(f,g)\mapsto B(f,g)$ is a product and in the product estimate one can choose the same constant as for $B$.

If $s\in(0,1), 1\leq p<\infty$ and $\varepsilon>0$ is such that $s+\varepsilon<1$ holds, then $B$ induces a product on $C^{s+\varepsilon}(\overline{I};X)\times W^s_p(I;Y)\rightarrow W^s_p(I;Z)$ and the constant in the estimate can be chosen independent of $|I|$ if $|I|\leq T$ for a fixed $T>0$.
\end{Lemma}
\textit{Proof.} First one can show that $B$ induces a product on $C^0(\overline{I};X)\times L^p(I;X)\rightarrow L^p(I;X)$ and in the estimate the same constant can be chosen. Then the first part follows using a null addition. Now let $f\in C^{s+\varepsilon}(\overline{I};X)$ and $g\in W^s_p(I;Y)$. Then a null addition yields
\[
[B(f,g)]^p_{W^s_p(I;Z)}\leq C_0^p([f]^p_{C^{s+\varepsilon}(\overline{I};X)}\int_I\int_I \|g(x)\|^p_{Y}\frac{dx\,dy}{|x-y|^{1-\varepsilon p}}+\|f\|^p_{C^0(\overline{I};X)}[g]^p_{W^s_p(I;Y)}).
\] 
By Tonelli and $I-x\subseteq[-|I|,|I|]$ the double integral can be estimated by $\|g\|_{L^p(I;Y)}^p \frac{2|I|^{\varepsilon p}}{\varepsilon p}$.\hfill$\square$

The following lemma yields embeddings for suitable intersection spaces:
\begin{Theorem}\label{th_intungl_einb_L_B}
Let $X_0,X_1,X$ be Banach spaces over $\K=\R$ or $\C$ with $X_1\hookrightarrow X\hookrightarrow X_0$ and $I\subseteq\R$ be an open interval. Assume that for a $\theta\in(0,1)$ the following interpolation inequality holds: $\|x\|_X\leq C_0\|x\|_{X_0}^{1-\theta}\|x\|_{X_1}^\theta$ for all $x\in X_0\cap X_1$. Then
\begin{enumerate}
\item If $1\leq p,p_0,p_1\leq\infty$ with $\frac{1}{p}=\frac{1-\theta}{p_0}+\frac{\theta}{p_1}$, then $L^{p_0}(I;X_0)\cap L^{p_1}(I;X_1)\hookrightarrow L^p(I;X)$.
\item If $1\leq p<\infty$, then $L^p(I;X_1)\cap W^1_p(I;X_0)\hookrightarrow B^{1-\theta}_{p,\infty}(I;X)$.
\end{enumerate}
In both cases analogous interpolation inequalities hold with $C_0$ and $3C_0$, respectively.
\end{Theorem}
\textit{Proof.} The first part follows from Hölder's inequality. Now let $f\in L^p(I;X_1)\cap W^1_p(I;X_0)$ and $0<h\leq 1$ be arbitrary. It remains to estimate $[f]_{B^{1-\theta}_{p,\infty}(I;X)}$. The first part yields
\[
\left\|\frac{f(.+h)-f}{h^{1-\theta}}\right\|_{L^p(I_h;X)}\leq
C_0\left\| \frac{f(.+h)-f}{h}\right\|^{1-\theta}_{L^p(I_h;X_0)}2^\theta\|f\|^\theta_{L^p(I;X_1)}.
\]
Using $f(x+h)-f(x)=h \int_0^1 f'(x+th)\,dt$, Hölder's inequality and Tonelli's theorem we get $\|f(.+h)-f\|^p_{L^p(I_h;X_0)}\leq|h|^p\|f'\|_{L^p(I;X_0)}^p$. Thus the claim follows.\hfill$\square$

The following theorem characterizes the trace of functions in $\E_1(T)$ at $t=0$.
\begin{Theorem}\label{th_einb_aus_trace_method}
Let $X_0,X_1$ Banach spaces over $\K=\R$ or $\C$ with $X_1\hookrightarrow X_0$ and let $1<p<\infty$. Then for $\E_1(\R_+):=L^p(\R_+;X_1)\cap W^1_p(0,\infty;X_0)$ it holds that
\[
X_\gamma:=(X_0,X_1)_{1-\frac{1}{p},p}=\{u(0):u\in\E_1(\R_+)\}
\] 
and an equivalent norm on $X_\gamma$ is given by the induced quotient norm. Moreover we have $\E_1(\R_+)\hookrightarrow C_b^0(\R_+;X_\gamma)$ and $\E_1(T)\hookrightarrow C^0([0,T];X_\gamma)$ for $0<T<\infty$ where the embedding constant is bounded independent of $T>0$ if we add $\|u(0)\|_{X_\gamma}$ in the $\E_1(T)$-norm. In particular this also holds for $u(0)=0$.
\end{Theorem}
\textit{Proof.} Up to $\E_1(\R_+)\hookrightarrow C_b^0(\R_+;X_\gamma)$ this directly follows from the trace method, cf. Lunardi \cite{Lunardi_Interpolation}, Section 1.2. The embedding $\E_1(T)\hookrightarrow C^0([0,T];X_\gamma)$ can be shown by applying the one on $\R_+$ to suitable extensions. E.g. one can first use reflection at $T$ and then multiply by a cutoff function. It remains to show the property of the embedding constant. In case $u(0)=0$ we extend $u$ as before and for $u(0)\neq 0$ we subtract an $\tilde{u}\in\E_1(\R_+)$ with $\tilde{u}(0)=u(0)$ such that $\|\tilde{u}\|_{\E_1(\R_+)}\leq 2\|u(0)\|_{1-\frac{1}{p},p}^{\tr}$ holds and apply the first part to $u-\tilde{u}$.\hfill$\square$

By interpolation we get a whole scale of embeddings from Theorem \ref{th_einb_aus_trace_method}:
\begin{Theorem}\label{th_einb_E_1T_in_hölder}
Let $X_0,X_1$ Banach spaces over $\K=\R$ or $\C$ with $X_1\hookrightarrow X_0$. Moreover, let $1<p<\infty$ and $T>0$. Then for $0<\theta<1-\frac{1}{p}$ and $1\leq q\leq\infty$ it holds
\[
\E_1(T)\hookrightarrow C^{1-\theta-\frac{1}{p}}([0,T];(X_0,X_1)_{\theta,q})
\]
and the embedding constant is bounded independent of $T>0$ if we add $\|u(0)\|_{X_\gamma}$ in the $\E_1(T)$-norm, where $\E_1(T)$ and $X_\gamma$ are as in Theorem \ref{th_einb_aus_trace_method}.
\end{Theorem}
\textit{Proof.} It is well-known that $W^1_p(0,T;X_0)\hookrightarrow C^{0,1-\frac{1}{p}}([0,T];X_0)$ and the semi-norm is estimated by $\|u'\|_{L^p(0,T;X_0)}$ for all $u\in W^1_p(0,T;X_0)$. Now one interpolates with the embedding from Theorem \ref{th_einb_aus_trace_method} and uses the reiteration theorem in Lunardi \cite{Lunardi_Interpolation}, Corollary 1.24 and the interpolation inequality in \cite{Lunardi_Interpolation}, Corollary 1.7.\hfill$\square$

As another consequence of Theorem \ref{th_einb_aus_trace_method} we obtain
\begin{Lemma}\label{th_einb_E_1T_mit_W2p_W1p0}
Let $\Omega\subseteq\R^n, n\geq 2$ be a bounded $C^\infty$-domain and $\Sigma\subseteq\Omega$ be a connected, compact and smooth hypersurface that separates $\Omega$ in two disjoint, connected domains $\Omega^\pm$ with $\Sigma=\partial\Omega^+$. Moreover, let $0<T<\infty, 2\leq p<\infty$ and $\Omega_0:=\Omega^+\cup\Omega^-$. Then
\[
\E_1(T):=L^p(0,T;X_1)\cap W^1_p(0,T;L^p(\Omega))\hookrightarrow C^0([0,T];Y_\gamma)
\]
with $X_1:=W^2_p(\Omega_0)\cap W^1_{p,0}(\Omega)$ and $Y_\gamma:=\{v\in L^p(\Omega):v|_{\Omega^\pm}\in W^{2-\frac{2}{p}}_p(\Omega^\pm)\}\cap W^1_{p,0}(\Omega)$. The embedding constant is bounded independent of $T>0$ if we add $\|u(0)\|_{Y_\gamma}$ in the norm.
\end{Lemma}
\textit{Proof.} We define $X_\gamma:=(L^p(\Omega),X_1)_{1-\frac{1}{p},p}$ and $\tilde{X}_\gamma:=(L^p(\Omega),W^2_p(\Omega_0))_{1-\frac{1}{p},p}$. Then Theorem \ref{th_einb_aus_trace_method} implies
$\E_1(T)\hookrightarrow C^0([0,T];X_\gamma)$ and $\tilde{\E}_1(T):=L^p(0,T;W^2_p(\Omega_0))\cap W^1_p(0,T;L^p(\Omega))$ embeds to $C^0([0,T];\tilde{X}_\gamma)$. The embedding constants are bounded in a similar way as in the theorem. E.g. from the $K$-method, cf. Lunardi \cite{Lunardi_Interpolation}, Section 1.1 we obtain
\[
\tilde{X}_\gamma=\{v\in L^p(\Omega):v|_{\Omega^\pm}\in W^{2-\frac{2}{p}}_p(\Omega^\pm)\}
\]
with equivalent norms. Hence it remains to show that for the continuous representative of $u\in\E_1(T)$ additionally $u(t)\in W^1_{p,0}(\Omega)$ for all $t\in[0,T]$ holds. Therefore we use the characterization in Theorem \ref{th_einb_aus_trace_method} to conclude $X_\gamma\hookrightarrow \tilde{X}_\gamma$. Here $X_1$ is dense in $X_\gamma$ by Lunardi \cite{Lunardi_Interpolation}, Proposition 1.17. For $v\in X_1$ it holds $\lsprung v\rsprung=\tr_\Sigma(v|_{\Omega^-})-\tr_{\Sigma}(v|_{\Omega^+})=0$ and $\tr_{\partial\Omega}\,v=0$. By density this carries over to all $v\in X_\gamma$ because of $X_\gamma\hookrightarrow\tilde{X}_\gamma\hookrightarrow W^1_p(\Omega_0)$ for $2\leq p<\infty$ and the continuity of the respective boundary trace operators. Using the definition of weak derivative and partial integration as well as the trace characterization of $W^1_{p,0}(\Omega)$ we get $X_\gamma\subseteq\tilde{X}_\gamma\cap W^1_{p,0}(\Omega)=Y_\gamma$.\hfill$\square$

Finally, we need some properties of certain types of intersection spaces:
\begin{Lemma}\label{th_W_C0}
Let $\Omega\subseteq \R^n$ a bounded, connected domain, $0<T\leq T_0$ and $1\leq p<\infty$. Then $
W_T:=C^0([0,T];C^0(\overline{\Omega}))\cap W^1_p(0,T;L^p(\Omega))$ is an algebra with pointwise multiplication and it holds $\|fg\|_{W_T}\leq C\|f\|_{W_T}\|g\|_{W_T}$ with $C>0$ independent of $T,T_0$ and $p$.

Moreover, for $f\in W_T$ with $f(t)(x)\geq c_0>0$ for all $(x,t)\in\overline{\Omega}\times [0,T]$ also $1/f\in W_T$ and we have $\partial_t(1/f)=-\partial_tf/f^2$. If additionally $\|f\|_{W_T}\leq R$, then $\|\frac{1}{f}\|_{W_T}\leq c_R$ with $c_R>0$ independent of $0<T\leq T_0$.
\end{Lemma}
\textit{Proof.} Using convolution one can show that $C^\infty([0,T];C^0(\overline{\Omega}))$ is dense in $W_T$. With this one can directly prove the first part and that a product rule holds. The second part follows by density utilizing Lemma \ref{th_hölder_ableiten} and that the evaluation at any point in $\overline{\Omega}$ gives a bounded, linear functional on $C^0(\overline{\Omega})$.\hfill$\square$\\

\section{Maximal regularity for a two-phase Stokes system}\label{sec_lin_stokes}
Let $2\leq q<3$, $\Omega\subseteq\R^n, n\geq 2$ be a bounded $C^\infty$-domain and $\Sigma\subseteq\Omega$ a compact, smooth and connected hypersurface that divides $\Omega$ into two disjoint and connected domains $\Omega^\pm$ with $\Sigma=\partial\Omega^+$ and outer unit normal $\nu_\Sigma$. We define $\Omega_0:=\Omega^+\cup\Omega^-$ and consider the following two-phase Stokes system
\begin{alignat}{2}\label{eq_stokes1}
\partial_tu-\mu^\pm\Delta u+\nabla \tilde{q}&=f &\qquad &\text{ in }\Omega^\pm\times(0,T),\\
\diverg\,u&=g &\qquad &\text{ in }\Omega^\pm\times(0,T),\label{eq_stokes2}\\
\lsprung u\rsprung&=0 &\qquad &\text{ on }\Sigma\times(0,T),\label{eq_stokes3}\\
\lsprung T(u,\tilde{q})\rsprung\,\nu_\Sigma&=a &\qquad &\text{ on }\Sigma\times(0,T),\label{eq_stokes4}\\
u|_{\partial\Omega}&=0 &\qquad &\text{ on }\partial\Omega\times(0,T),\label{eq_stokes5}\\
u|_{t=0}&=u_0 &\qquad &\text{ in }\Omega,\label{eq_stokes6}
\end{alignat}
where $\mu^\pm>0$ and $T(u,\tilde{q}):=2\mu^\pm\textup{sym}(\nabla u)-\tilde{q}I$ in $\Omega^\pm$. With the aid of Shimizu \cite{Shimizu} we show an optimal existence result in an $L^q$-setting: For $(u,\tilde{q})$ we introduce the space $Y_T:=Y_T^1\times Y_T^2$ where
\begin{align*}
Y_T^1&:= W^1_q(0,T;L^q(\Omega))^n\cap L^q(0,T;W^2_q(\Omega_0)\cap W^1_{q,0}(\Omega))^n,\\
Y_T^2&:=\{\tilde{q}\in L^q(0,T;W^1_{q,(0)}(\Omega_0)):\lsprung \tilde{q} \rsprung\in W^{1-\frac{1}{q},\frac{1}{2}(1-\frac{1}{q})}_q(\Sigma_T)\}.
\end{align*}
Here we require mean value $0$ for $\tilde{q}$ in order to get uniqueness later. This implies necessary conditions for the data $(f,g,a,u_0)$ if \eqref{eq_stokes1}-\eqref{eq_stokes6} hold:
\begin{Lemma}\label{th_stokes_lin_notw}
Let $2\leq q<3,T>0$ and $(u,\tilde{q})\in Y_T$. Then $(f,g,a,u_0)$ defined by \eqref{eq_stokes1}-\eqref{eq_stokes6} is contained in $\tilde{Z}_T:=\tilde{Z}_T^1\times \tilde{Z}_T^2\times (\tilde{Z}_T^3)^n\times Y_\gamma$ with $\tilde{Z}_T^1:=L^q(0,T;L^q(\Omega))^n$ and
\begin{align*}
\tilde{Z}_T^2&:=\{ g\in L^q(0,T;W^1_{q,(0)}(\Omega_0))\cap W^1_q(0,T;W^{-1}_{q,(0)}(\Omega)):
\tr_\Sigma (g|_{\Omega^+})\in \tilde{Z}_T^3\},\\
\tilde{Z}_T^3&:= W_q^{1-\frac{1}{q},\frac{1}{2}(1-\frac{1}{q})}(\Sigma_T),\quad Y_\gamma:=\{v\in L^q(\Omega):v|_{\Omega^\pm}\in W^{2-\frac{2}{q}}_q(\Omega^\pm)\}^n\cap W^1_{q,0}(\Omega)^n.
\end{align*} 
Additionally, the compatibility condition $\diverg\,u_0=g|_{t=0}$ in $W^{-1}_{q,(0)}(\Omega)$ is valid and we have $\|(f,g,a,u_0)\|_{\tilde{Z}_T}\leq C(\|(u,\tilde{q})\|_{Y_T}+\|u_0\|_{Y_\gamma})$ with $C>0$ independent of $(u,\tilde{q})$ and $T>0$.
\end{Lemma}
\begin{Remark}\label{th_bem_stokes_lin}\upshape
The spaces are equipped with the natural norms. In particular we include the term $\|\lsprung\tilde{q}\rsprung\|_{\tilde{Z}_T^3}$ in the $Y_T^2$-norm and $\|\tr_\Sigma(g|_{\Omega^+})\|_{\tilde{Z}_T^3}$ in the $\tilde{Z}_T^2$-norm. One can directly verify that all spaces are Banach spaces.
\end{Remark}
\textit{Proof of Lemma \ref{th_stokes_lin_notw}.} The assertion for $f$ and $L^q$-in-time properties of $g$ can be directly shown. Lemma \ref{th_einb_E_1T_mit_W2p_W1p0} implies $u_0\in Y_\gamma$. Weak differentiability of $g$ follows from
\[
\int_\Omega g(t)\phi\,dx=-\int_\Omega u(t)\cdot\nabla\phi\,dx
\]
for all $\phi\in W^1_{q'}(\Omega)$ and almost all $t\in(0,T)$. The latter identity also implies the compatibility condition. Results on mixed order Sobolev spaces yield together with $\lsprung\tilde{q}\rsprung\in\tilde{Z}^3_T$ that $a, \tr_\Sigma(g|_{\Omega^+})\in\tilde{Z}_T^3$ and the uniform estimate in the lemma can be shown by extending functions in time with suitable operators, cf. proof of Lemma \ref{th_G}, below \eqref{eq_bew_G_4}.\hfill$\square$

Our maximal regularity result for \eqref{eq_stokes1}-\eqref{eq_stokes6} is
\begin{Theorem}\label{th_stokes_lin}
Let $2\leq q<3$ and $0<T\leq T_0$. Let $(f,g,a,u_0)$ be contained in the space $Z_T:=\{(f,g,a,u_0)\in\tilde{Z}_T:\diverg\,u_0=g|_{t=0}\text{ in }W^{-1}_{q,(0)}(\Omega)\}$. Then the Stokes system \eqref{eq_stokes1}-\eqref{eq_stokes6} has exactly one solution $(u,\tilde{q})$ in $Y_T$. Moreover, there is a $C>0$ independent of $(f,g,a,u_0)$ and $0<T\leq T_0$ such that $\|(u,\tilde{q})\|_{Y_T}\leq C\|(f,g,a,u_0)\|_{\tilde{Z}_T}$.
\end{Theorem}
\begin{Remark}\label{th_bem_stokes_lin2}\upshape
\begin{enumerate}
\item Further necessary conditions are $\tr_{\partial\Omega}\,g\in W^{1-\frac{1}{q},\frac{1}{2}(1-\frac{1}{q})}_q(\partial\Omega_T)$ and $\tr_\Sigma (g|_{\Omega^-})\in\tilde{Z}_T^3$. But these follow indirectly from the proof of Theorem \ref{th_stokes_lin}.
\item For simplicity we restricted to $2\leq q<3$. In principle also for other scales of $q$ statements are possible. But then e.g. other necessary conditions can arise. For instance, if $q>3$ another compatibility condition appears involving the jump of the stress tensor, cf. Shimizu \cite{Shimizu}.
\end{enumerate}
\end{Remark}
\textit{Proof of Theorem \ref{th_stokes_lin}.} The case $(a,g,u_0)=0$ follows from the announced result by Shimizu \cite{Shimizu}, Theorem 2. The main step for the proof is the investigation of model problems. For $\R^n$ with $\R^{n-1}$ as interface cf. Shibata and Shimizu \cite{Shibata_Shimizu_Model_Prob}. For simplicity we do not give a proof of Theorem 2 in \cite{Shimizu}. But we reduce the general case to the above one similar to A. and Wilke \cite{Abels_Mullins_Sekerka}, proof of Theorem A.1.

First of all, we reduce to the case $(g,u_0)|_{\Omega^+}=0$. Therefore we consider
\begin{alignat}{2}\label{eq_bew_mr_stklin1}
\diverg\,v^+&=g|_{\Omega^+}&\qquad&\text{ in }\Omega^+\times(0,T),\\
v^+|_{t=0}&=u_0|_{\Omega^+}&\qquad&\text{ in }\Omega^+\label{eq_bew_mr_stklin2}
\end{alignat}
and apply Theorem 1.1 in A. \cite{Abels_nonstat_Stokes} (for $\Gamma_1=\emptyset,\Gamma_2=\Sigma$ and $f,a=0$ there). To this end we have to show the properties of $g|_{\Omega^+}$ and $u_0|_{\Omega^+}$. First it holds
\[
g|_{\Omega^+}\in L^q(0,T;W^1_q(\Omega^+))\quad\text{ and }\quad\tr_\Sigma(g|_{\Omega^+})\in W^{1-\frac{1}{q},\frac{1}{2}(1-\frac{1}{q})}_{q}(\Sigma_T)
\] 
by definition of $\tilde{Z}_T^2$. Moreover, for $w\in L^q_{(0)}(\Omega)$ and $\phi\in W^1_{q',0}(\Omega^+)$ 
\[
\langle w|_{\Omega^+},\phi\rangle_{W^{-1}_{q}(\Omega^+),W^1_{q',0}(\Omega^+)}=\int_{\Omega^+}w\phi\,dx=\int_\Omega w\,\tilde{e}_0\phi\,dx=\langle w,\tilde{e}_0\phi\rangle_{W^{-1}_{q,(0)}(\Omega),W^1_{q',(0)}(\Omega)},
\]
where $e_0\phi\in W^1_{q',0}(\Omega)$ is the extension by zero of $\phi$ and $\tilde{e}_0\phi:=e_0\phi-\frac{1}{|\Omega|}\int_{\Omega}(e_0\phi)(y)\,dy$. In particular this is also valid for $g(t)$ instead of $w$ for almost all $t\in(0,T)$. Since
\[
W^{-1}_{q,(0)}(\Omega)\rightarrow W^{-1}_{q}(\Omega^+):w\mapsto [\phi\mapsto\langle w,\tilde{e}_0\phi\rangle]
\] 
is bounded and linear, we obtain $g|_{\Omega^+}\in W^1_q(0,T;W^{-1}_{q}(\Omega^+))$, the norm is estimated by $C\|g\|_{\tilde{Z}_T^2}$ and for all $\phi\in W^1_{q',0}(\Omega^+)$ we have
\[
\langle (g|_{\Omega^+})|_{t=0},\phi\rangle_{W^{-1}_{q}(\Omega^+),W^1_{q',0}(\Omega^+)}=\langle g|_{t=0},\tilde{e}_0\phi\rangle_{W^{-1}_{q,(0)}(\Omega),W^1_{q',(0)}(\Omega)}.
\]
This yields $\diverg(u_0|_{\Omega^+})=(g|_{\Omega^+})|_{t=0}$ in $W^{-1}_{q}(\Omega^+)$ because of $u_0\in Y_\gamma$ and the compatibility condition in $Z_T$. Now Theorem 1.1 in A. \cite{Abels_nonstat_Stokes} yields a solution $v^+\in W^{2,1}_q(\Omega^+_T)^n$ of \eqref{eq_bew_mr_stklin1}-\eqref{eq_bew_mr_stklin2} with
\[
\|v^+\|_{W^{2,1}_{q}(\Omega^+_T)^n}\leq C_{T_0}(\|g\|_{\tilde{Z}_T^2}+\|u_0\|_{Y_\gamma}),
\]
where $C_{T_0}>0$ is independent of $0<T\leq T_0$. We extend $v^+$ to $\tilde{v}^+\in Y^1_T$ using a suitable extension operator. Because of Lemma \ref{th_stokes_lin_notw} we have reduced to the case $(g,u_0)|_{\Omega^+}=0$.

Now we want to trace the latter one back to the case $(g,u_0)=0$. Therefore we look at
\begin{alignat}{2}\label{eq_bew_mr_stklin3}
\diverg\,v^-&=g|_{\Omega^-}&\qquad&\text{ in }\Omega^-\times(0,T),\\
v^-|_{\partial\Omega^-}&=0&\qquad&\text{ on }\partial\Omega^-\times(0,T),\label{eq_bew_mr_stklin4}\\
v^-|_{t=0}&=u_0|_{\Omega^-}&\qquad&\text{ in }\Omega^-\label{eq_bew_mr_stklin5}
\end{alignat}
and again apply Theorem 1.1 in A. \cite{Abels_nonstat_Stokes} (for $\Gamma_1=\partial\Omega^-,\Gamma_2=\emptyset$ and $f,a=0$ there). Since $g|_{\Omega^+}=0$ it holds $
g|_{\Omega^-}\in L^q(0,T;W^1_{q,(0)}(\Omega^-))$. Furthermore, for $w\in L^q_{(0)}(\Omega)$ with $w|_{\Omega^+}=0$ and all $\phi\in W^1_{q'}(\Omega^-)$ we have
\[
\langle w|_{\Omega^-},\phi\rangle_{W^1_{q'}(\Omega^-)',W^1_{q'}(\Omega^-)}=\int_{\Omega^-}w\phi\,dx=\int_{\Omega}w\,\tilde{E}\phi\,dx=\langle w,\tilde{E}\phi\rangle_{W^{-1}_{q,(0)}(\Omega),W^1_{q',(0)}(\Omega)},
\]
where $E:W^1_{q'}(\Omega^-)\rightarrow W^1_{q'}(\Omega)$ is an extension operator and $\tilde{E}\phi:=E\phi-\frac{1}{|\Omega|}\int_\Omega(E\phi)(y)\,dy$. Since 
\[
W^{-1}_{q,(0)}(\Omega)\rightarrow W^1_{q'}(\Omega^-)':w\mapsto[\phi\mapsto\langle w,\tilde{E}\phi\rangle]
\]
is bounded and linear, we obtain $g|_{\Omega^-}\in W^1_q(0,T;W^1_{q'}(\Omega^-)')$, the norm is estimated by $C\|g\|_{\tilde{Z}_T^2}$ and for all $\phi\in W^1_{q'}(\Omega^-)$ it holds that
\[
\langle(g|_{\Omega^-})|_{t=0},\phi\rangle_{W^1_{q'}(\Omega^-)',W^1_{q'}(\Omega^-)}=\langle g|_{t=0},\tilde{E}\phi\rangle_{W^{-1}_{q,(0)}(\Omega),W^1_{q',(0)}(\Omega)}.
\]
Since $u_0|_{\Omega^+}=0$ and $u_0\in Y_\gamma$ it follows that $\diverg(u_0|_{\Omega^-})=(g|_{\Omega^-})|_{t=0}$ in $W^1_{q'}(\Omega^-)'$ with the compatibility condition in $Z_T$. By Theorem 1.1 in A. \cite{Abels_nonstat_Stokes} there exists a solution $v^-\in W^{2,1}_{q}(\Omega^-_T)^n$ of \eqref{eq_bew_mr_stklin3}-\eqref{eq_bew_mr_stklin5} such that
\[
\|v^-\|_{W^{2,1}_{q}(\Omega^-_T)^n}\leq C_{T_0}(\|g\|_{\tilde{Z}_T^2}+\|u_0\|_{Y_\gamma}),
\] 
where $C_{T_0}>0$ is independent of $0<T\leq T_0$. Because of $v^-|_{\partial\Omega^-}=0$, we can extend $v^-$ by zero to $e_0(v^-)\in Y^1_T$. Lemma \ref{th_stokes_lin_notw} implies that we have reduced to the case $(g,u_0)=0$.

To trace this one back to the case $(g,a,u_0)=0$ we define for any $a:\Sigma\rightarrow\R^n$ the normal and tangential component $a_\nu:=a\cdot\nu_\Sigma\,\nu_\Sigma$ and $a_\tau:=a-a_\nu$, respectively. First we reduce to the case $(g,a_{\tau},u_0)=0$. As in A. \cite{Abels_nonstat_Stokes}, proof of Lemma 2.5, 2. there is an $A\in W^{2,1}_{q}(\Omega^+_T)^n$ such that $\|A\|_{W^{2,1}_{q}(\Omega^+_T)^n}\leq C\|a\|_{(\tilde{Z}_T^3)^n}$ as well as
\[
A|_{t=0}=0,\quad A|_{\Sigma}=0,\quad (2\mu^+\textup{sym}(\nabla A)|_\Sigma \nu_\Sigma)_\tau=a_\tau\quad\text{ and }\quad (\diverg\,A)|_\Sigma=0.
\]
But we also want $A$ to be divergence free. To this end one shows as in the proof of Lemma \ref{th_stokes_lin_notw} using $A|_\Sigma=0$ that $\diverg\,A\in L^q(0,T;W^1_{q,0}(\Omega^+))\cap W^1_q(0,T;W^1_{q'}(\Omega^+)')$ and the norm is estimated by $C\|A\|_{W^{2,1}_q(\Omega^+)^n}$. Now we apply the Bogovski\u{\i}-operator $B$ to $\diverg\,A$, cf. Geissert, Heck and Hieber \cite{Bogovskii}, Theorem 2.5, and obtain
\[
B(\diverg\,A)\in L^q(0,T;W^2_{q,0}(\Omega^+))^n\cap W^1_q(0,T;L^q(\Omega^+))^n\quad\text{ with }\quad\diverg(B(\diverg\,A))=\diverg\,A,
\]
since $\diverg\,A$ has mean value $0$. Additionally, the norm of $B(\diverg A)$ is estimated by $C\|A\|_{W^{2,1}_q(\Omega^+)^n}$. Therefore the extension by zero to $\Omega$ of $\tilde{A}:=A-B(\diverg\,A)$ has the desired properties and we reduced to the case $(g,a_\tau,u_0)=0$.

Finally, to reduce the latter one to the case $(g,a,u_0)=0$, we subtract a suitable extension of $a\cdot\nu_\Sigma\in\tilde{Z}_T^3$ from the pressure: By the trace theorem there is a
\[
p^+\in L^q(0,T;W^1_q(\Omega^+))\quad\text{ with }\|p^+\|_{L^q(0,T;W^1_q(\Omega^+))}\leq C\|a\|_{(\tilde{Z}_T^3)^n}\quad\text{ and }\quad p^+|_\Sigma=a\cdot\nu_\Sigma.
\]
We extend $p^+$ by $0$ to $\Omega$ and subtract the mean value which does not alter the jump. Hence we obtain a $\tilde{p}^+\in Y_T^2$ with $\|\tilde{p}^+\|_{Y_T^2}\leq C\|a\|_{(\tilde{Z}_T^3)^n}$ and $-\lsprung\tilde{p}^+\rsprung\,\nu_\Sigma=a_\nu$.

Altogether we reduced the general case to the case $(g,a,u_0)=0$.\hfill$\square$\\

\setcounter{secnumdepth}{0}\setcounter{tocdepth}{0} 


\end{document}